\documentclass[11pt,leqno]{article}

\usepackage{amstex,amscd,amssymb,theorem,epsfig}

\long\def\comment#1\endcomment{}

\makeatletter
\begingroup
\gdef\th@dotted{\normalfont\itshape
  \def\@begintheorem##1##2{%
        \item[\hskip\labelsep \theorem@headerfont ##1\ ##2.]}%
\def\@opargbegintheorem##1##2##3{%
   \item[\hskip\labelsep \theorem@headerfont ##1\ ##2\ (##3).]}}
\endgroup
\makeatother

\theoremstyle{dotted}

\newtheorem{thm}{Theorem}[section]
\newtheorem{lemma}[thm]{Lemma}

\newtheorem{prop}[thm]{Proposition}

\makeatletter
\begingroup
\gdef\th@upshape{\normalfont
  \def\@begintheorem##1##2{%
        \item[\hskip\labelsep \theorem@headerfont ##1\ ##2.]}%
\def\@opargbegintheorem##1##2##3{%
   \item[\hskip\labelsep \theorem@headerfont ##1\ ##2\ (##3).]}}
\endgroup
\makeatother

\theoremstyle{upshape}

\newtheorem{defn}[thm]{Definition}
\newtheorem{rem}[thm]{Remark}

\makeatletter
\renewcommand{\subsection}{\@startsection{subsection}{2}{0pt}{-3ex
plus -1ex minus -0.2ex}{-2mm plus -0pt minus
-2pt}{\normalfont\bfseries}} 
\makeatother

\newcounter{figuredammit}  

\makeatletter
\@addtoreset{equation}{section}
\@addtoreset{figure}{section}
\@addtoreset{figuredammit}{section}
\makeatother

\newcommand{\proof}[1][Proof.]{\smallskip\noindent{\em #1}}
\def\endproof{\hfill\ensuremath{\square}\par\medskip}

\def\eqref#1{\thetag{\ref{#1}}}

\let\latexref=\ref
\def\ref#1{{\normalfont{\latexref{#1}}}}

\newcommand{\wt}{\widetilde}
\newcommand{\6}{\partial}

\newcommand{\cntrct}
{\hspace{2pt}\raisebox{1pt}{\text{$\lrcorner$}}\hspace{2pt}}

\newcommand{\hdot}{{\:\raisebox{3pt}{\text{\circle*{1.5}}}}}
\newcommand{\h}{{\mathbb H}}
\newcommand{\C}{{\mathbb C}}
\newcommand{\R}{{\mathbb R}}
\newcommand{\Z}{{\mathbb Z}}

\newcommand{\LL}{{\cal L}}
\newcommand{\CC}{{\cal C}}

\newcommand{\T}{{\cal T}}
\newcommand{\E}{{\cal E}}
\newcommand{\F}{{\cal F}}
\newcommand{\B}{{\cal B}}

\newcommand{\conj}{\overline{\phantom{a}}}

\newcommand{\wgt}{\operatorname{\sf wt}}
\newcommand{\Alt}{\operatorname{\sf Alt}}

\newcommand{\Ker}{\operatorname{Ker}}
\newcommand{\id}{\operatorname{\sf id}} 

\begin{document}

\title{A canonical hyperk\"ahler metric on the total space of a
cotangent bundle} 

\author{D. Kaledin}

\maketitle

\begin{abstract}
A canonical hyperk\"ahler metric on the total space $T^*M$ of a
cotangent bundle to a complex manifold $M$ has been constructed
recently by the author in \cite{K}. This paper presents the results
of \cite{K} in a streamlined and simplified form. The only new
result is an explicit formula obtained for the case when $M$ is an
Hermitian symmetric space.
\end{abstract}

\section*{Introduction.}

Constructing a hyperk\"ahler metric on the total space $T^*M$ of the
cotangent bundle to a K\"ahler manifold $M$ is an old problem,
dating back to the very first examples of hyperk\"ahler metrics
given by E. Calabi in \cite{C}. Since then, many people have
obtained a lot of important results valid for manifolds $M$ in this
or that particular class (see, for example, the papers \cite{Sw},
\cite{BG}, \cite{K2}, \cite{K2}, \cite{Nak}). Finally, the general
problem has been more or less solved a couple of years ago,
independently by B. Feix \cite{F} and by the author \cite{K}.

The metrics constructed in \cite{F} and \cite{K} are the same. In
fact, this metric satisfies an additional condition which makes it
essentially unique -- which justifies the use of the term
``canonical metric''.  But the approaches in \cite{F} and \cite{K}
are very different. Feix's method is very geometric in nature; it is
based on a direct description of the associated twistor space. The
approach in \cite{K} is much farther from geometric
intuition. However, it seems to be more likely to lead to explicit
formulas.

Unfortunately, the paper \cite{K} is 100 pages long, and it is not
written very well. Some of the proofs are not at all easy to
understand and to check. The exposition is sometimes canonical to
the point of obscurity.

Recently the author has been invited to give a talk on the results
of \cite{K} at the Second Quaternionic Meeting in Rome. It was a
good opportunity to revisit the subject and to streamline and
simplify some of the proofs. This paper, written for the Proceedings
volume of the Rome conference, is an attempt to present the results
of \cite{K} in a concrete and readable form. The exposition is
parallel to \cite{K} but the paper is mostly independent.

Compared to \cite{K}, the emphasis in the present paper has been
shifted from canonical but abstract constructions to things more
explicit and down-to-earth. The number of definitions is reduced to
the necessary minimum. I have also tried to give a concrete
geometric interpretation to everything that admits such an
interpretation. In a sense, the exposition as compared to \cite{K}
intentionally goes to the other extreme. Thus the present paper is
not so much a replacement for \cite{K} but rather a companion paper
-- the same story told in a different way.

I should note that the canonical hyperk\"ahler metric on $T^*M$ has
one important defect -- namely, it is defined only in an open
neighborhood $U \subset T^*M$ of the zero section $M \subset
T^*M$. This raises a very interesting and difficult ``convergence
problem''. One would like to describe the maximal open subset $U
\subset T^*M$ where the canonical metric is defined, and to say when
$U$ is the whole total space $T^*M$. Unfortunately, very little is
known about this. In fact, in the present paper we even restrict
ourselves to giving the formal germ of the canonical metric near $M
\subset T^*M$. The fact that this formal germ converges to an actual
metric at least on an open subset $U \subset T^*M$ is proved in the
last Section of \cite{K}. The proof is long and tedious but
completely straightforward. Since I do not know how to improve it, I
have decided to omit it altogether to save space.

This reader will find the precise statements of all the results in
Section~\ref{deffs}. The last part of that Section contains a brief
description of the rest of the paper, and indicates the parallel
places in \cite{K}, the differences in notation and terminology and
so on. The only thing in this paper which is completely new is the
last Section~\ref{symmm}. It contains an explicit formula for the
canonical metric (or rather, for the canonical hypercomplex
structure) in the case when $M$ is an Hermitian symmetric space. The
formula is similar to the general formula for symmetric $M$ obtained
by O. Biquard and P. Gauduchon in \cite{BG}.

\bigskip

\noindent
{\bf Acknowledgments.} I would like to thank the organizers of the
Rome Quaternionic Meeting for inviting me to this very interesting
conference and for giving me an opportunity to present the results
of \cite{K}. Part of the present work was done during my visit to
Ecole Polytechnique in Paris during the Autumn of 1999. I have
benefited a lot from the hospitality and the stimulating atmosphere
of this institution. I would like to thank P. Gauduchon for inviting
me to Paris and for encouraging me to write up a streamlined version
of \cite{K}. The present paper owes a lot to discussions with
O. Biquard and P. Gauduchon during my visit. In particular, the last
Section is an attempt to compare \cite{K} with the results in their
beautiful paper \cite{BG}. I am also grateful to M. Verbitsky, who
read the first draft of the manuscript and suggested several
improvements.

\section{Statements and definitions.}\label{deffs}

To save space, we will assume some familiarity with hyperk\"ahler
and hypercomplex geometry. We only give a brief reminder. The reader
will find excellent expositions of the subject in \cite{B},
\cite{HKLR}, \cite{S1}, \cite{S2}.

Let $\h$ be the algebra of quaternions. A smooth manifold $X$ is
called {\em almost quaternionic} if it is equipped with a smooth
action of the algebra $\h$ on the tangent bundle $TX$. Equivalently,
one can consider a smooth action on the cotangent bundle
$\Lambda^1(M)$. To fix terminology, we will assume that $\h$ acts on
$\Lambda^1(M)$ on the left.

An almost quaternionic manifold $X$ is called {\em hypercomplex} if
it admits a torsion-free connection preserving the $\h$-module
structure on $\Lambda^1(M)$. Such a connection exists is unique. It
is called the {\em Obata connection} of the hypercomplex manifold
$X$.

A Riemannian almost quaternionic manifold $X$ is called {\em
hyper-hermitian} if the Riemannian pairing satisfies
$$
\left(h\alpha_1,\alpha_2\right) =
\left(\alpha_1,\overline{h}\alpha_2\right), \qquad\qquad
\alpha_1,\alpha_2 \in \Lambda^1(M), h \in \h,
$$
where $\overline{h}$ is the quaternion conjugate to $h$. A
hyper-hermitian almost quaternionic manifold $X$ is called {\em
hyperk\"ahler} if the $\h$-action is parallel with respect to the
Levi-Civita connection $\nabla_{LC}$. In other words, $X$ must be
hypercomplex, and the Obata connection must be $\nabla_{LC}$.

Let $X$ be a almost quaternionic manifold. Every embedding $\C
\hookrightarrow \h$ from the field of complex numbers to the algebra
of the quaternions induces an almost complex structure on $X$. If
$X$ is hypercomplex, then all these induced almost complex
structures are integrable. If $X$ is also hyperk\"ahler, then all
these complex structures are K\"ahler with respect to the metric.

Throughout this paper it will be convenient to choose an embedding
$I:\C \to \h$ and an additional element $j \in \h$ such that 
\begin{align*}
j^2 &= -1\\
j \cdot I(z) &= I(\overline{z}) \cdot j, \qquad z \in \C
\end{align*}
With these choices, every left $\h$-module $V_\R$ defines a complex
vector space $V = V_I$ and a map $j:V \to \overline{V}$ which
satisfies
\begin{equation}\label{j.square}
j \circ \overline{j} = -\id
\end{equation}
We will call $V_I$ the {\em main complex structure} on the real
vector space $V_\R$. Conversely, every pair $\langle V, j\rangle$ of
a complex vector space $V$ and a map $j:V \to \overline{V}$ which
satisfies \eqref{j.square} defines an $\h$-module structure on the
real vector space $V_\R$ underlying $V$.

The map $j:V \to \overline{V}$ can be considered as an automorphism
$J:V_\R \to V_\R$ of the underlying real vector space. This map
induces a second complex structure on $V_\R$. We will call it the
{\em complementary} complex structure and denote the resulting
complex vector space by $V_J$

Applying this to vector bundles, we see that an almost quaternionic
manifold $X$ is the same as an almost complex manifold $X$ equipped
with a smooth complex bundle map
\begin{equation}\label{j.def}
j: \T(X) \to \overline{\T}(X)
\end{equation}
from the tangent bundle $\T(X)$ to its complex-conjugate bundle
$\overline{\T}(X)$ which satisfies \eqref{j.square}.

It would be very convenient to have some way to know whether an
almost quaternionic manifold $X$ is hypercomplex or hyperk\"ahler
without working explicitly with torsion-free connections. For
hypercomplex manifolds, the integrability condition is very
simple. An almost quaternionic manifold $X$ is hypercomplex if and
only if both the main complex structure on $X$ and the complementary
almost complex structure $X_J$ are integrable.

The simplest way to determine whether a hyper-hermitian almost
quaternionic manifold $X$ is hyperk\"ahler is to consider the
complex-bilinear form $\Omega$ on the tangent bundle $\T(X)$ given
by
\begin{equation}\label{omega}
\Omega(\xi_1,\xi_2) = h(\xi_1,j(\xi_2)), \qquad \xi_1,\xi_2 \in
\T(X),
\end{equation}
where $h$ is the Hermitian metric on $X$. Then \eqref{j.square}
insures that the form $\Omega$ is skew-symmetric. The manifold $X$
is hyperk\"ahler if and only if both the $(2,0)$-form $\Omega$ and
the K\"ahler $\omega$ are closed. Alternatively, one can define a
$(2,0)$-form $\Omega_J$ for the complementary almost complex
structure $X_J$ instead of the $(2,0)$-form $\Omega$ for the main
almost complex structure. Then $X$ is hyperk\"ahler if and only if
both $\Omega$ and $\Omega_J$ are closed. Moreover, it suffices to
require that these $(2,0)$-forms are holomorphic, each in its
respective almost complex structure on $X$ (``holomorphic'' here
means that $d\Omega$ is a form of type $(3,0)$). Finally, if we
already know that the manifold $X$ is hypercomplex, then it suffices
to require that only one of the forms $\Omega$, $\Omega_J$ is a
holomorphic $2$-form.

If $X$ is a K\"ahler manifold equipped with a closed $(2,0)$-form
$\Omega$, one can {\em define} a map $j:\T(X) \to \overline{\T}(X)$
by \eqref{omega}. Then $X$ is hyperk\"ahler if and only if this map
$j$ satisfies \eqref{j.square}.

Consider now the case when $X = T^*M$ is the total space of the
cotangent bundle to a K\"ahler manifold $M$ with the k\"ahler metric
$h$. Then $X$ carries a canonical holomorphic $2$-form
$\Omega$. Moreover, the unitary group $U(1)$ acts on $X$ by
dilatations along the fibers of the projection $X \to M$, so that we have
\begin{equation}\label{omega.weight.one}
z^* \Omega = z\Omega
\end{equation}
for every $z \in U(1) \subset \C$. Using these data, we can
formulate the main result of \cite{K} as follows.

\begin{thm}\label{hk}
There exists a unique, up to fiber-wise automorphisms of $X/M$,
$U(1)$-invariant K\"ahler metric $h$ on $X =T^*M$, defined in the
formal neighborhood of the zero section $M \subset T^*M = X$, such
that
\begin{enumerate}
\item $h$ restricts to the given K\"ahler metric on the zero section
$M \subset X$, and
\item the pair $\langle \Omega, h \rangle$ defines a hyperk\"ahler
structure on $X$ near $M \subset X$.
\end{enumerate}
Moreover, if the K\"ahler metric $h$ on $M$ is real-analytic, then
the formal canonical metric on $X$ comes from a real-analytic
hyperk\"ahler metric defined in an open neighborhood of $M \subset
X$.
\end{thm}

The $U(1)$-action on $T^*M$ is very important for this theorem. In
fact, $T^*M$ with this action belongs to a general class of
hyperk\"ahler manifolds equipped with a $U(1)$-action, introduced in
\cite{H}, \cite{HKLR}.

\begin{defn}\label{comp}
We will say that a holomorphic $U(1)$-action on a
hy\-per\-k\"ah\-ler manifold $X$ is {\em compatible with the
hyperk\"ahler structure} if and only if
\begin{enumerate}
\item the metric $h$ on $X$ is $U(1)$-invariant,
\item the holomorphic $2$-form $\Omega$ satisfies
\eqref{omega.weight.one}, 
\item the map $j:\T(X) \to \overline{\T}(X)$ satisfies
\begin{equation}\label{j.weight.one}
j(z^*\xi) = z \cdot z^*(j(\xi)), \qquad \xi \in \T(X), z \in U(1)
\subset \C.
\end{equation}
\end{enumerate}
\end{defn}

It is easy to check using \eqref{omega} that every two of these
conditions imply the third.

Theorem~\ref{hk} can generalized to the following statement,
somewhat analogous to the Darboux-Weinstein Theorem in symplectic
geometry. 

\begin{thm}\label{darboux}
Let $X$ be a hyperk\"ahler manifold equipped with a regular
compatible holomorphic $U(1)$-action. Then there exists an open
neighborhood $U \subset X$ of the $U(1)$-fixed point subset $M =
X^{U(1)} \subset X$ and a canonical embedding $\LL:U \to T^*M$ such
that the hyperk\"ahler structure on $U$ is induced by means of the
map $\LL$ from the canonical hyperk\"ahler structure on $T^*M$.
\end{thm}

Here {\em regular} is a certain condition on the $U(1)$-action near
the fixed points subset $X^{U(1)} \subset X$ which is formulated
precisely in Definition~\ref{reg} (roughly speaking, weights of the
action on the tangent space $T_mX$ at every point $m \in X^{U(1)}
\subset X$ should be $0$ and $1$). The map $\LL$ will be called the
{\em normalization map}. Note that this Theorem allows one to
reformulate Theorem~\ref{hk} so that the metric is indeed unique --
not just unique up to a fiber-wise automorphism of $X/M$. To fix the
metric, it suffices to require that the associated normalization map
$\LL:X \to X = \overline{T}M$ is identical. Metrics with this
property will be called {\em normalized}.

There is also a form of Theorem~\ref{hk} for hypercomplex manifolds
(and it is this form which is the most important for \cite{K} -- all
the other statements are obtained as its corollaries). To formulate
it, we note that out of the three condition of
Definition~\ref{comp}, the third one makes sense for almost
quaternionic (in particular, hypercomplex) manifolds. We will say
that a holomorphic $U(1)$-action on an almost quaternionic manifold
$X$ is {\em compatible with the quaternionic action} if
Definition~\ref{comp}~\thetag{iii} is satisfied.

Let $\overline{T}M$ be the total space of the bundle
$\overline{\T}(M)$ complex-conjugate to the tangent bundle of the
manifold $M$. Then $\overline{T}M$ is a smooth manifold, and we have
the canonical projection $\rho:\overline{T}M \to M$ and the zero
section $i:M \to \overline{T}M$. The group $U(1)$ acts on
$\overline{T}M$ by dilatations along the fibers of the projection
$\rho$. Moreover, for any compatible hypercomplex structure on the
$U(1)$-manifold $\overline{T}M$ the corresponding Obata connection
$\nabla_O$ induces a torsion-free connection $\nabla$ on $M$ by the
following rule
$$
\nabla(\alpha) = i^*(\nabla_O\rho^*\alpha), \qquad \alpha \in
\Lambda^1(M).
$$
The hypercomplex version of Theorem~\ref{hk} is the following. 

\begin{thm}\label{hc}
Let $M$ be a complex manifold $M$ equipped with a holomorphic connection
$\nabla$ on the tangent bundle $\T(M)$ such that
\begin{enumerate}
\item $\nabla$ has no torsion, and
\item the curvature of the connection $\nabla$ is of type $(1,1)$.
\end{enumerate}
Let $X = \overline{T}M$ be the total space of the complex-conjugate
to the tangent bundle $TM$. Let the group $U(1)$ act on $X$ by
dilatations along the fibers of the projection $\rho:X \to M$.

Then there exists a unique, up to a fiber-wise automorphism of $X/M$,
hypercomplex structure on $X$, defined in the formal neighborhood of
the zero section $M \subset X$, such that the embedding $i:M
\hookrightarrow X$ and the projection $\rho:X \to M$ are holomorphic
and the Obata connection on $X$ induces the given connection
$\nabla$ on $M$.

Moreover, if the connection $\nabla$ is real-analytic, then the
hypercomplex structure on $X$ is real-analytic in an open
neighborhood $U \subset X$ of $M \subset X$.
\end{thm}

Note that {\em a priori} there is no natural complex structure on
the space $X = \overline{T}M$ (we write $\overline{T}$ instead of
$T$ just to indicate the correct $U(1)$-action). Therefore
Theorem~\ref{hc} in fact claims two things: firstly, there exists an
integrable almost complex structure on $X$, and secondly, there
exists a map $j:\T(X) \to \overline{\T}(X)$ which extends it to a
hypercomplex structure.

Connections $\nabla$ that satisfy the conditions of this Theorem
were called {\em k\"ahlerian} in \cite{K} (see \cite[8.1.2]{K}). I
would like to thank D. Joyce for attracting my attention to his
paper \cite{J}, where he uses the same class of connections to
construct {\em commuting} almost complex structures on the total
space $TM$ of the tangent bundle to a complex manifold $M$. Joyce
calls these connection {\em K\"ahler-flat}.

Theorem~\ref{hc} also admits a generalized Darboux-like version in
the spirit of Theorem~\ref{darboux} (see \cite[Proposition 4.1]{K});
we do not formulate it here to save space.

\bigskip

We will now give a brief outline of the remaining part of the
paper. In Section~\ref{nrm} we consider an arbitrary
$U(1)$-equivariant hypercomplex manifold $X$ and construct the
normalization map $\LL:X \to \overline{T}M$, thus proving
Theorem~\ref{darboux}. This corresponds to \cite[Section 4]{K}. What
we call {\em normalization} here was called {\em linearization} in
\cite{K}; {\em normalized} corresponds to {\em linear}. The
terminology of \cite{K} has been changed because it was misleading:
connections on the fibration $\overline{T}M \to M$ that were called
linear in \cite{K} are not linear in the usual sense of the word.

Section~\ref{hb} introduces $\R$-Hodge structures and the so-called
Hodge bundles (Definition~\ref{hb.def}) which are the basis of our
approach to $U(1)$-equivariant quaternionic manifolds $X$. This
corresponds to \cite[Sections 2,3]{K}. Proposition~\ref{dol} is a
version of \cite[Proposition 3.1]{K}.

In Section~\ref{hcon} we turn to the case $X = \overline{T}M$ and
introduce the so-called Hodge connections
(Definition~\ref{hcon.def}). This corresponds to \cite[Section
5]{K}. The proofs have been considerably shortened. There are also
some new facts on the relation between our formalism and the objects
one usually associates with hypercomplex manifolds. In particular,
Lemma~\ref{j.matrix} is new.

Section~\ref{wa} and Section~\ref{pf} introduce the Weil algebra
$\B^\hdot(M)$ of a complex manifold $M$ (Definition~\ref{wa.def})
and then use it to prove Theorem~\ref{hc}. Section~\ref{wa} contains
the preliminaries; Section~\ref{pf} gives the proof itself. This
material corresponds to \cite[Sections 6-8]{K}. The approach has
been changed in the following way. Keeping track of various gradings
and bigradings and on the Weil algebra presents considerable
difficulties: when the proof of Theorem~\ref{hc} is written down,
the number of indices becomes overwhelming. In \cite{K} we have
tried to handle this by an auxiliary technical device called {\em
the total Weil algebra} (\cite[7.2.4]{K}). It was quite a natural
thing to do from the conceptual point of view. Unfortunately, the
proof became more abstract than one would like. Here we have opted
for the direct approach. To make things comprehensible, we rely on
pictures (Figure~\ref{fig.odd}, Figure~\ref{fig.even}) which
graphically represent the Hodge diamonds of the relevant pieces of
the Weil algebra $\B^\hdot(M)$.

Finally, Section~\ref{metrics} deals with things hyperk\"ahler: we
deduce Theorem~\ref{hk} from Theorem~\ref{hc}. This corresponds to
\cite[Section 9]{K}. We believe that the exposition has also been
simplified, and the proofs are easier to check.

The last Section~\ref{symmm} of this paper is new. We try to
illustrate our constructions by a concrete example of an Hermitian
symmetric space $M$. We obtain a formula similar to \cite{BG}. The
last section of \cite{K} contains the proof of convergence of our
formal series in the case when the K\"ahler manifold $M$ is
real-analytic. In this paper, this proof is entirely omitted.

\section{Normalization.}\label{nrm}

Of all the statement formulated in the last Section, the most
straightforward one is the Darboux-like Theorem~\ref{darboux} and
its hypercomplex version. In this Section, we explain how to
construct the normalization map $\LL$. Most of the proofs are
replaced with references to \cite{K}.

We begin with some generalities. Assume that the group $U(1)$ acts
smoothly on a smooth manifold $X$. For any point $x \in X$ fixed by
the action, we have an action of $U(1)$ on the tangent space
$T_xX$. Equivalently, we have the weight decomposition 
$$
T_xX \otimes_\R \C = \bigoplus_k (T_xX)^k,
$$
where $z \in U(1) \subset \C$ acts on $(T_xX)^k$ by multiplication
by $z^k$. We will say that the fixed point $x$ is {\em regular} if
the only non-trivial weight components $(T_xX)^k$ correspond to
weights $k = 0,1$. The subset $X^{U(1)} \subset X$ of fixed points
is a disjoint union of smooth submanifolds of different dimensions.
Regular fixed points form a connected component of this set. Denote
this component by $M \subset X$.

Let $\phi$ be the differential of the $U(1)$-action -- that is, the
vector field on $X$ which gives the infinitesimal action of the
generator $\frac{\6}{\6\theta}$ of the Lie algebra of the group
$U(1)$. Assume further that $X$ is a complex manifold and that
$U(1)$ preserves the complex structure. Say that a point $x \in X$
is {\em stable} if for any $t \in \R$, $t \geq 0$ there exists
$\exp(\sqrt{-1}t\phi)x$, and moreover, the limit
$$
x_0 = \lim_{t \to +\infty}\exp(\sqrt{-1}t\phi)x
$$
also exist. When the limit point $x_0$ does exist, it is obviously
fixed by $U(1)$. Say that a stable point $x \in X$ is {\em regular
stable} if the limit point is a regular fixed point, $x_\infty \in M
\subset X$. Regular stable points form an open subset $X^{reg}
\subset X$.

\begin{defn}\label{reg}
A complex $U(1)$-manifold $X$ is called is {\em regular} if every
point $x \in X$ is regular stable, $X^{reg} = X$.
\end{defn}

For example, the total space of an arbitrary complex vector bundle
on an arbitrary complex manifold is regular (if $U(1)$ acts by
dilatations along the fibers). The submanifold of regular fixed
points in this example is the zero section.

When the $U(1)$-manifold $X$ is hypercomplex, we will say that it is
regular if it is regular in the main complex structure. In this
case, the subset $M \subset X$ of regular fixed points is a complex
submanifold. Setting
$$
\rho(x) = x_0 = \lim_{t \to +\infty}\exp(\sqrt{-1}t\phi)x
$$
defines a $U(1)$-invariant projection $\rho:X \to M$.

\begin{lemma}[{{\cite[4.2.1-3]{K}}}]
The projection $\rho:X \to M$ is a holomorphic submersion.\endproof
\end{lemma}

We can now define the normalization map. Consider the exact sequence
$$
\begin{CD}
0 @>>> \T(X/M) @>>> \T(X) @>>> \rho^*\T(M) @>{d\rho}>> 0
\end{CD}
$$
of tangent bundles associated to the submersion $\rho:X \to M$. The
differential $\phi$ of the $U(1)$-action is a vertical holomorphic
vector field on $X$, $\phi \in \T(X/M)$. Applying the operator
$j:\T(X) \to \overline{\T}(X)$ to $\phi$ gives a section of the
bundle $\overline{\T}(X)$. We can project this section to obtain a
section
$$
d\rho(j(\phi)) \in \rho^*\overline{\T}(M)
$$
of the pullback bundle $\rho^*\overline{\T}(M)$. But such a section
tautologically defines a map $\LL:X \to \overline{T}M$ from $X$ to
the total space of the complex bundle $\overline{T}$ on the manifold
$M$. Since $\phi$ is $U(1)$-invariant, and $j$ is of weight $1$ with
respect to the $U(1)$-action, the section $d\rho(j(\phi))$ is also
of weight one. This means that the associated map
$$
\LL:X \to \overline{T}M
$$
is $U(1)$-equivariant. We will call it the {\em normalization map}
for the regular hypercomplex $U(1)$-manifold $X$.

\begin{lemma}[{{\cite[Proposition 4.1]{K}}}]
The normalization map $\LL:X \to \overline{T}M$ is an open
embedding. \endproof
\end{lemma}

This Lemma essentially reduces Theorem~\ref{darboux} to
Theorem~\ref{hk}.

A particular case occurs when $X = \overline{T}M$ is itself the
total space of the complex-conjugate to the tangent bundle on a
complex manifold $M$ -- or, more generally, an open
$U(1)$-invariant neighborhood $U \subset \overline{T}M$ of the zero
section $M \subset \overline{T}M$. As noted above, in this case the
zero section $M \subset \overline{T}(M)$ coincides with the subset
of fixed points. Therefore the normalization map $\LL:U \to
\overline{T}M$ is an open embedding from $U$ into $\overline{T}M$,
possibly different from the given one.

\begin{defn}
The hypercomplex structure on $U \subset \overline{T}M$ is called
{\em normalized} if the normalization map $\LL:U \to \overline{T}M$
coincides with the given embedding.
\end{defn}

(In particular, when $U = \overline{T}M$ is the whole total space,
the normalization map must be identical.)

To prove Theorem~\ref{hk} and Theorem~\ref{hc}, it is sufficient to
be able to classify all normalized hypercomplex structures on the
$U(1)$-manifold $\overline{T}(M)$ and germs of such structures near
the zero section $M \subset \overline{T}M$. It will be convenient to
slightly rewrite the normalization condition. Namely, the identity
map $\id:\overline{T}M \to \overline{T}M$ defines a section on $X =
\overline{T}(M)$ of the pullback bundle $\rho^*\overline{\T}(M)$. We
will denote this section by $\tau$ and call it the {\em tautological
section}. Then a hypercomplex structure on $U \subset X$ is
normalized if and only if we have
\begin{equation}\label{norm}
j(\phi) = \tau \in \rho^*\overline{\T}(M).
\end{equation}

\section{Hodge bundles.}\label{hb}

The first step in the proof of Theorem~\ref{hc} is to give a
workable description of hypercomplex structures on the total space
$X = \overline{T}M$. For this we use the language of $\R$-Hodge
structures. 

Recall that an {\em $\R$-Hodge structure $V$ of weight $k$} is by
definition a real vector space $V_\R$ equipped with a grading 
\begin{equation}\label{grr}
V = V_\R \otimes_\R \C = \bigoplus_p V^{p,k-p}
\end{equation}
such that
\begin{equation}\label{cnj}
\overline{V^{p,q}} = V^{q,p}, \qquad p, q \in \Z, p+q = k.
\end{equation}
Equivalently, instead of the grading \eqref{grr} one can consider a
$U(1)$ action on $V$ defined by
$$
z \cdot v = z^pV, \qquad v \in V^{p,q} \subset V, z \in U(1) \subset
\C.
$$
Then \eqref{cnj} becomes
$$
\overline{z \cdot v} = z^k z \cdot \overline{v}, \qquad v \in V, z
\in U(1) \subset \C.
$$
When the weight $k$ is equal to $1$, this equation on the complex
conjugation map becomes precisely \eqref{j.weight.one}. The
difference between the complex conjugation map and the map $j$ used
to define quaternionic structures is that the first one is an
involution, $\overline{\overline{v}} = v$, while for the second one
we have $j(j(v)) = -v$. Nevertheless, we will exploit the similarity
between them to describe quaternionic actions via Hodge
structures. To do this, we use the following trick. Let $V$ be an
$\R$-Hodge structure of weight $1$, and consider the map
$$
\iota:V \to V
$$
given by the action of $-1 \in U(1) \subset \C$ -- in other words,
let 
$$
\iota(v) = (-1)^p, \qquad v \in V^{p,1-p} \subset V.
$$
Then the map
\begin{equation}\label{j.conj}
j = \iota \circ \conj:V \to \overline{V}
\end{equation}
still satisfies \eqref{j.weight.one}, and \eqref{j.square} also
holds. This turns $V$ into a left $\h$-module. Conversely, every
left $\h$ module $\langle V, j\rangle$ equipped with a $U(1)$-action
on $V$ such that $j$ satisfies \eqref{j.weight.one} defines an
$\R$-Hodge structure of weight $1$.

To use this for a description of hypercomplex structures on
manifolds, we introduce the following.

\begin{defn}\label{hb.def}
Let $X$ be a smooth manifold equipped with an action of the group
$U(1)$, and let $\iota:X \to X$ be the action of the element $-1 \in
U(1) \subset \C$. A {\em Hodge bundle $\E$ of weight $k$} on $X$ is
by definition a $U(1)$-equivariant complex vector bundle $\E$
equipped with a $U(1)$-equivariant isomorphism
$$
\conj:\E \to \overline{\E}(k)
$$
such that $\conj \circ \conj = \id$.
\end{defn}

Here $\overline{\E}(k)$ is the bundle complex conjugate to $\E$,
whose $U(1)$-equivariant structure is twisted by tensoring with the
$1$-dimensional representation $\C(k)$ of the group $U(1)$ of weight
$k$, 
$$
z \cdot x = z^kx, \qquad x \in \C(k), z \in U(1) \subset \C.
$$
When the $U(1)$-action on the manifold $X$ is trivial, a weight-$k$
Hodge bundle $\E$ on $X$ is just the bundle of $\R$-Hodge structures
of weight $k$ in the obvious sense. In particular, if $X$ is an
almost complex manifold, then the bundle $\Lambda^k(X)$ of all
complex-valued $k$-forms on $X$ is a Hodge bundle of weight
$k$.

When the $U(1)$-action on $X$ is no longer trivial, every bundle
$\E$ of $\R$-Hodge structures on $X$ still defines a Hodge
bundle. Thus $\Lambda^k(X)$ is still a weight-$k$ Hodge bundle. But
this Hodge bundle structure is not interesting, since it does not
take into account the natural $U(1)$-action on
$\Lambda^k(X)$. Assuming that the $U(1)$-action preserves the almost
complex structure on $X$, we define instead a Hodge bundle structure
on $\Lambda^1(X,\C)$ by keeping the usual complex conjugation map
and twisting the $U(1)$-action so that
$$
\Lambda^1(X,\C) \cong \Lambda^{1,0}(X)(1) \oplus \Lambda^{0,1}(X)
$$
as a $U(1)$-equivariant vector bundle. It is easy to check that this
indeed defines on $\Lambda^1(X,\C)$ a Hodge bundle structure of
weight $1$.

Assume now that the almost complex manifold $X$ is equipped with an
almost quaternionic structure which is compatible with the
$U(1)$-action. Then the complex vector bundle $\Lambda^{0,1}(X)$ of
$(0,1)$-forms on $X$ already has a structure of a Hodge bundle of
weight $1$. The $U(1)$-action for this structure is the natural one,
and the complex conjugation map is induced by the map $j:\T(X) \to
\overline{\T}X$ via \eqref{j.conj}. This Hodge bundle structure
completely determines the quaternionic action. More precisely, for
every smooth $U(1)$-manifold $X$, every Hodge bundle $\E$ of weight
$1$ whose underlying real vector bundle $\E_\R$ is identified with
the cotangent bundle $\Lambda^1(X,\R)$ comes from a unique
compatible almost quaternionic structure $X$.

The natural embedding $\Lambda^{0,1}(X) \subset \Lambda^1(X,\C)$ is
not a map of Hodge bundles -- it is $U(1)$-equivariant, but it
obviously does not commute with the complex conjugation map. This
can be corrected. To do this, one has to look at the picture in a
different way (which will turn out to be very useful). Return for a
moment to linear algebra. Let $V_\R$ be a left $\h$-module, and let
$V$, $V_J$ be the complex vector spaces obtained from $V_\R$ by the
main and the complementary complex structures. Consider the complex
vector space $V_\R \otimes_\R \C$. This vector space does not depend
on the $\h$-action on $V_\R$. Given an $\h$-action, we have the main
and the complementary complex structure operators $I =
I(\sqrt{-1}):V_\R \to V_\R$ and $J:V_\R \to V_\R$ and the associated
eigenspace decompositions
\begin{align}
V_\R \otimes_\R \C &= V \oplus \overline{V},\label{main.spl}\\
V_\R \otimes_\R \C &= V_J \oplus \overline{V_J}\label{comp.spl}
\end{align}
Since the operators $I$ and $J$ anti-commute, these decompositions
are distinct: we have
$$
V \cap V_J = V \cap \overline{V_J} = \overline{V} \cap V_J =
\overline{V} \cap \overline{V_J} = 0.
$$
In particular, the composition
\begin{equation}\label{H}
H:\overline{V} \to V_\R \otimes \C \to \overline{V_J}
\end{equation}
of the canonical embedding in \eqref{main.spl} and the canonical
projection in \eqref{comp.spl} is an isomorphism. We will call it
the {\em canonical isomorphism between the main and the
complementary complex structures}. On the level of the real vector
space $V_\R$, the map $H$ is induced by a non-trivial automorphism
$H:V_\R \to V_\R$ (in fact it is the action of the element
$I(\sqrt{-1}) + j \in \h$). Conjugation with this map interchanges
the operators $I$ and $J$.

Return now to the case of an almost quaternionic manifold $X$. Then
we claim that the complementary almost complex structure operator
$J:\Lambda^1(X,\C) \to \Lambda^1(X,\C)$ is a map of Hodge
bundles. Indeed, it commutes with the complex conjugation map by
definition. Therefore it suffices to show that it is
$U(1)$-equivariant on $\Lambda^{0,1}(X) \subset \Lambda^1(X,\C)$.
But for every $v \in \Lambda^{1,0}(X)$ we have $J(v) = j(v) \in
\Lambda^{0,1}(X)$, and the map $j$ is of weight $1$. Thus the
operator $J$ is indeed $U(1)$-equivariant (recall the twisting of
the $U(1)$-action in the definition of the Hodge bundle structure on
$\Lambda^1(X,\C)$).

Since the endomorphism $J:\Lambda^1(X,\C) \to \Lambda^1(X,\C)$ is a
map of Hodge bundles, its eigenbundles 
$$
\Lambda^{1,0}_J(X),\Lambda^{0,1}_J(X) \subset \Lambda^1(X,\C)
$$
are Hodge subbundles. Therefore we obtain a canonical weight-$1$
Hodge bundle structure on the bundle $\Lambda^{0,1}_J(X)$ of
$(0,1)$-forms for the {\em complementary} almost complex structure
on $X$. 

This Hodge bundle is not a new one. Indeed, the canonical
isomorphism 
$$
H:\Lambda^{0,1}(X) \to \Lambda^{0,1}_J(X)
$$
defined in \eqref{H} is $U(1)$-equivariant -- it is obtained as a
composition of $U(1)$-equivariant maps. Moreover, it is very easy to
check that $H$ commutes with the complex conjugation. Thus
$\Lambda^{0,1}(X) \cong \Lambda^{0,1}_J(X)$ as Hodge bundles. But
the projections
\begin{align}
\Lambda^1(X,\C) &\to \Lambda^{0,1}(X),\\
\Lambda^1(X,\C) &\to \Lambda^{0,1}_J(X) \cong
\Lambda^{0,1}(X),\label{prj} 
\end{align}
are different. Only the second one is a Hodge bundle map.

All this linear algebra is somewhat tautological, but it becomes
useful when we consider the integrability conditions on an almost
quaternionic quaternionic structure. The real advantage of
$\R$-Hodge structures over $\h$-modules is the presence of higher
weights. Namely, the category of $\h$-modules admits no natural
tensor product. On the other hand, the category of $\R$-Hodge
structures and the category of Hodge bundles are obviously tensor
categories. Thus, for example, the weight-$1$ Hodge bundle structure
on the cotangent bundle $\Lambda^1(X,\C)$ induces a weight-$k$ Hodge
bundle structure on the bundle $\Lambda^k(X,\C)$ of $k$-forms.

To make use of these higher-weight Hodge bundles, we need a
convenient notion of maps between Hodge bundles of different
weights.

\begin{defn}
A bundle map (or, more generally, a differential operator) $f:\E \to
\F$ between Hodge bundles $\E$, $\F$ of weights $m$, $n$ is called
{\em weakly Hodge} if it commutes with the complex conjugation map
and admits a decomposition
\begin{equation}\label{Htype}
f = \sum_{0 \leq p \leq n-m} f^p,
\end{equation}
where $f^p:\E \to \F$ is of weight $p$ with respect to the
$U(1)$-action -- in other words, $f^p$ is $U(1)$-equivariant when
considered as a map
$$
f^p:\E \to \F(p).
$$
\end{defn}

We see that non-trivial weakly Hodge maps between Hodge bundles
$\E$, $\F$ exist only when their weights satisfy $\wgt\F \geq
\wgt\E$.

When the $U(1)$-action on the manifold $X$ is trivial, the Hodge
bundles $\E$ and $\F$ come from bundles of $\R$-Hodge structures on
$X$, and the decomposition $f = \sum_p f^p$ of a weakly Hodge map
$f:\E \to \F$ is simply the Hodge type decomposition,
$$
f^p = f^{p,q}, \qquad p + q = \wgt\F - \wgt \E.
$$
If the $U(1)$-action is not trivial, but preserves an almost complex
structure on $X$, the de Rham differential $d = \6 + \bar\6 =
d^{1,0} + d^{0,1}:\Lambda^0(X) \to \Lambda^1(X)$ is weakly Hodge. If
the almost complex structure is integrable, then the same is true
for the de Rham differential $d:\Lambda^k(X,\C) \to
\Lambda^{k+1}(X,\C)$ for every $k \geq 0$.

When the $U(1)$-manifold is almost quaternionic, we have a Hodge
bundle structure of weight $1$ on the bundle
$\Lambda^{0,1}_J(X)$. Then the Dolbeault differential
$$
D = \bar\6_J:\Lambda^0(X,\C) \to \Lambda^{0,1}_J(X)
$$
is weakly Hodge. Indeed, it is the composition of the weakly Hodge
de Rham differential and the projection \eqref{prj} which is a Hodge
bundle map. Let
\begin{equation}\label{d0}
D = D^0 + D^1
\end{equation}
be the decomposition \eqref{Htype} for the weakly Hodge map
$D$. Looking at the definition of the canonical isomorphism
$H:\Lambda^{0,1}(X) \to \Lambda^{0,1}_J(X)$, we see that the
component 
$$
D^0:\Lambda^0(X,\C) \to \Lambda^{0,1}_J(X)
$$ 
in the decomposition \eqref{d0} coincides with the Dolbeault
differential
$$
\6:\Lambda^0(X,\C) \to \Lambda^{0,1}(X) \cong \Lambda^{0,1}_J(X)
$$
for the main complex structure.

Assume now that the almost quaternionic $U(1)$-manifold $X$ is
hypercomplex. The bundle $\Lambda^{0,k}_J(x)$ of $(0,k)$-forms on
$X_J$ is a Hodge bundle of weight $k$ for every $k \geq 0$, and we
have the Dolbeault differential
$$
D = \bar\6_J:\Lambda^{0,k}_J(x) \to \Lambda^{0,k+1}_J(x).
$$
Since the projections $\Lambda^k(X,\C) \to \Lambda^{0,k}_J(X)$ are
Hodge bundle maps for every $k \geq 0$, this Dolbeault differential
is weakly Hodge. It turns out that this is a sufficient
integrability condition for an almost quaternionic manifold equipped
with a compatible $U(1)$-action.

\begin{prop}\label{dol}
Let $X$ be an almost quaternionic manifold equipped with a
compatible $U(1)$-action. Assume that the Dolbeault differential
$$
D:\Lambda^0(X,\C) \to \Lambda^{0,1}_J(X)
$$
extends to a weakly Hodge derivation $D:\Lambda^{0,\hdot}_J(x) \to
\Lambda^{0,\hdot+1}_J(X)$ of the algebra $\Lambda^{0,\hdot}_J(X)$
satisfying $D \circ D = 0$. Then the manifold $X$ is hypercomplex.
\end{prop}

\proof{} It suffices to prove that both the main and the
complementary almost complex structures on $X$ are integrable. For
this, it is enough to prove that the Dolbeault differentials
\begin{align*}
\bar\6_J = D&:\Lambda^0(X,\C) \to \Lambda^{0,1}_J(X),\\
\bar\6 = D^0&:\Lambda^0(X,\C) \to \Lambda^{0,1}_J(X) \cong
\Lambda^{0,1}(X),
\end{align*}
extend to square-zero derivations of the exterior algebra
$\Lambda^{0,\hdot}_J(X)$. The differential $D$ extends by
assumption. To extend $D^0$, take the component
$D^0:\Lambda^{0,\hdot}_J(x) \to \Lambda^{0,\hdot+1}_J(X)$ of the
weakly Hodge map $D:\Lambda^{0,\hdot}_J(x) \to
\Lambda^{0,\hdot+1}_J(X)$. Then $D^0 \circ D^0$ is a component in
the decomposition \eqref{Htype} of the weakly Hodge map $D \circ
D:\Lambda^{0,\hdot}_J(x) \to \Lambda^{0,\hdot+2}_J(X)$. Since $D
\circ D = 0$, we also have $D^0 \circ D^0 = 0$.\endproof

We will now say a couple of words about hyperk\"ahler manifolds and
Hodge bundles. Let $X$ be an almost quaternionic
$U(1)$-manifold. Then every Riemannian metric on $X$ defines a
$(2,0)$-form $\Omega_J \in \Lambda^{2,0}(X)$. It turns out that if
the metric is hyper-hermitian and $U(1)$-invariant, then in
terminology of \cite{K}, the form $\Omega_J$ {\em of $H$-type
$(1,1)$}. This means the following. Consider the form $\Omega_J$ as
a bundle map
\begin{equation}\label{ht}
\Omega_{J}:\R(-1) \to \Lambda^{2,0}(X),
\end{equation}
where $\R(-1)$ is the trivial bundle on $X$ equipped with the
so-called Hodge-Tate $\R$-Hodge structure of weight $-1$ -- that is,
the complex conjugation map on $\R(-1)$ is minus the complex
conjugation map on $\R$, and the $U(1)$-equivariant structure is
twisted by $1$. The form $\Omega_J$ is said to be of $H$-type
$(1,1)$ if the map \eqref{ht} is a Hodge bundle map.

Conversely, every $(2,0)$-form $\Omega_J \in \Lambda^{2,0}(X)$ of
$H$-type $(1,1)$ on an almost quaternionic $U(1)$-manifold $X$ which
satisfies a positivity condition
\begin{equation}\label{posit}
\Omega(\xi_1,I(\xi_2)) > 0, \qquad\qquad \xi_1,\xi_2 \in T(X,\R)
\end{equation}
defines a $U(1)$-invariant hyper-hermitian metric on $X$. (See
\cite[1.5.4]{K}, but the proof is almost trivial.)

If $X$ is hypercomplex, then, as indicated in Section~\ref{deffs},
the metric corresponding to such a form $\Omega_J$ is hyperk\"ahler
if and only if the form $\Omega_J$ is holomorphic, $D\Omega_J = 0
\in \Lambda^{2,1}_J(X)$.

\begin{rem}
In fact, using the $U(1)$-action on $X$, one can even drop the
integrability condition. Indeed, if a form $\Omega_J$ of $H$-type
$(1,1)$ on an almost quaternionic $U(1)$-manifold $X$ satisfies
$D\Omega_J = 0$, then it also must satisfy
$$
D^0\Omega_J = D^1\Omega_J = 0.
$$
The canonical endomorphism $H:\Lambda^1(X,\R) \to \Lambda^1(X,\R)$,
being the conjugation with a quaternion, preserves up to a
coefficient the metric associated to $\Omega_J$ and interchanges the
almost complex structure operators $I$ and $J$. Therefore it sends
$\Omega_J$ to a form proportional to $\Omega$. Then $D^0\Omega_J=0$
implies that not only $\Omega_J$ is holomorphic, but that $\Omega$
is holomorphic as well. This proves that $X$ is hyperk\"ahler ({\em
a posteriori}, also hypercomplex). We will never need nor use this
argument. An interested reader will find details in \cite[3.3.4]{K}.
\end{rem}

\section{Hodge connections.}\label{hcon}

We will now restrict our attention to the case when the
$U(1)$-manifold $X$ is the total space $\overline{T}M$ of the
complex-conjugate to the tangent bundle of a complex manifold
$M$. In this case, Proposition~\ref{dol} is really useful, because
it turns out that the Hodge bundle algebra $\Lambda^{0,k}_J(X)$ does
not depend on an almost quaternionic structure on $X$. To see this,
denote by $\rho:X = \overline{T}M \to M$ the canonical projection,
and let
$$
\delta\rho:\rho^*\Lambda^1(M,\C) \to \Lambda^1(X,\C)
$$ 
be the codifferential of the map $\rho$. Then for every compatible
hypercomplex structure on $X$, we have canonical bundle maps
$$
\begin{CD}
\rho^*\Lambda^1(M,\C) @>{\delta\rho}>> \Lambda^1(X,\C) @>>>
\Lambda^{0,1}_J(X)
\end{CD}
$$
Assume that the manifold $X$ is a equipped with a hypercomplex
structure satisfying the conditions of Theorem~\ref{hc} -- namely,
assume that the projection $\rho:X \to M$ and the zero section $i:M
\to X$ are holomorphic maps. Then the codifferential $\delta\rho$ is
obviously a map of Hodge bundles. However, we also have the
following.

\begin{lemma}[{{\cite[5.1.9-10]{K}}}]\label{isom}
The composition map
$$
\rho^*\Lambda^1(M,\C) \to \Lambda^{0,1}_J(X)
$$
is an isomorphism of Hodge bundles in an open neighborhood of the
zero section $M \subset X$.\endproof
\end{lemma}

From this point on, it will be convenient to only consider germs of
hypercomplex structures defined near the zero section $M \subset
X$. In other words, we replace $X = \overline{T}M$ with an
unspecified and shrinkable $U(1)$-invariant open neighborhood of
the zero section. Since we are only interested in hypercomplex
structures on $X$ that satisfy the conditions of Theorem~\ref{hc},
Lemma~\ref{isom} shows that no matter what the particular
hypercomplex structure on $X$ is, we can {\em a priori} canonically
identify the Hodge bundle $\Lambda^{0,1}_J(X)$ with the pullback
bundle $\rho^*\Lambda^1(M,\C)$. The only thing that depends on the
hypercomplex structure is the derivation $D:\Lambda^0(X,\C) \to
\Lambda^{0,1}_J(X) \cong \rho^*\Lambda^1(M,\C)$.

To formalize the situation, we introduce the following.

\begin{defn}\mbox{}\label{hcon.def}
\begin{enumerate}
\item A {\em $\C$-valued connection} $\Theta$ on $X/M$ is a bundle map
$$
\Theta:\Lambda^1(X,\C) \to \rho^*\Lambda^1(M,\C)
$$
which splits the codifferential $\delta\rho:\rho^*\Lambda^1(M,\C)
\to \Lambda^1(X,\C)$ of the projection $\rho:X \to M$.

\item The {\em derivation $D:\Lambda^0(X,\C) \to
\rho^*\Lambda^1(M,\C)$ associated} to a $\C$-valued connection
$\Theta$ is the composition $D = \Theta \circ d$ of $\Theta$ with
the de Rham differential $d$.

\item A {\em Hodge connection} $\Theta$ on $X/M$ is a $\C$-valued
connection such that the associated derivation $D:\Lambda^0(X,\C)
\to \rho^*\Lambda^1(M,\C)$ is weakly Hodge.

\item A Hodge connection $\Theta$ on $X/M$ is called {\em flat} if
the associated derivation $D$ extends to a weakly Hodge derivation
$$
D:\rho^*\Lambda^\hdot(M,\C) \to \rho^*\Lambda^{\hdot+1}(M,\C)
$$
of the pullback of the de Rham algebra $\Lambda^\hdot(M,\C)$, and
the extended map $D$ satisfies $D \circ D = 0$.
\end{enumerate}
\end{defn}

Of course, a Hodge connection $\Theta$ is completely defined by the
associated derivation $D:\Lambda^0(X,\C) \to
\rho^*\Lambda^1(M,\C)$. Conversely, an arbitrary weakly Hodge
derivation
$$
D:\Lambda^0(X,\C) \to \rho^*\Lambda^1(M,\C)
$$ 
comes from a Hodge connection if and only if we have
$$
D\rho^*f = \rho^*df \in \rho^*\Lambda^1(M,\C)
$$
for every smooth function $f \in \Lambda^0(M,\C)$. Say that a
derivation 
$$
D:\Lambda^0(X,\C) \to \rho^*\Lambda^1(M,\C)
$$ 
is {\em holonomic} if the associated Hodge connection $\Theta$
induces an isomorphism
\begin{equation}\label{hl}
\Theta:\Lambda^1(X,\R) \to \rho^*\Lambda^1(M,\C)
\end{equation}
between the real cotangent bundle $\Lambda^1(X,\R)$ and the real
bundle underlying the complex vector bundle
$\rho^*\Lambda^1(M,\C)$. Then Proposition~\ref{dol} and
Lemma~\ref{isom} show that hypercomplex structures on $X$ satisfying
the conditions of Theorem~\ref{hc} are in one-to-one correspondence
with Hodge connections on $X/M$ whose associated derivations are
holonomic. Indeed, the isomorphism \eqref{hl} induces a Hodge bundle
structure of weight $1$ on the cotangent bundle $\Lambda^1(X,\R)$,
hence an almost quaternionic structure on $X$. Applying
Proposition~\ref{dol}, we see that flatness of the Hodge connection
is equivalent to the integrability of this almost quaternionic
structure. All derivations $D$ that we will work with will be
automatically holonomic -- this will turn out to be a consequence of
the normalization condition \eqref{norm} (see Lemma~\ref{D.norm}).

The name ``Hodge connection'' invokes the notion of a connection on
a smooth fibration. This is somewhat misleading. The problem is
that a Hodge connection $\Theta:\Lambda^1(X,\C) \to
\rho^*\Lambda^1(M,\C)$ is only defined over $\C$. So it has a real
part $\Theta_{Re}$ and an imaginary part $\Theta_{Im}$. The real
part
$$
\Theta_{Re}:\Lambda^1(X,\R) \to \rho^*\Lambda^1(M,\R)
$$
is indeed a connection on the fibration $\rho:X \to M$ in the usual
sense -- that is, it defines a smooth splitting 
$$
\Lambda^1(X,\R) \cong \rho^*\Lambda^1(M,\R) \oplus \Ker\Theta_{Re}
$$
of the real cotangent bundle $\Lambda^1(X,\R)$ into a horizontal and
a vertical part. The vertical part $\Ker\Theta_{Re}$ is canonically
isomorphic to the relative cotangent bundle $\Lambda^1(X/M,\R)$.

The imaginary part $\Theta_{Im}$, on the other hand, vanishes on the
subbundle $\rho^*\Lambda^1(M,\R) \subset \Lambda^1(X,\R)$ and
defines therefore a certain map
\begin{equation}\label{Rj.eq}
R_J:\Lambda^1(X/M,\R) \to \rho^*\Lambda^1(M,\R)
\end{equation}
from the relative cotangent bundle $\Lambda^1(X/M,\R)$ to the
pullback bundle $\rho^*\Lambda^1(M)$.

Since $X$ is an open subset in $\overline{T}M$, we can canonically
identify the bundle $\Lambda^1(X/M)$ with the pullback bundle
$\rho^*\Lambda^1(M)$. Under this identification, the map $R_J$
becomes an endomorphism of the bundle $\rho^*\Lambda^1(M)$.

Typically, when a Hodge connection $\Theta$ comes from a
hypercomplex structure on $X$, the associated real connection
$\Theta_{Re}$ on $X/M$ is {\em not} flat. It is only the sum
$$
\Theta = \Theta_{Re} + \sqrt{-1}\Theta_{Im}
$$
which is flat -- but it is no longer a real connection. This
situation is somewhat similar to what happens in C. Simpson's theory
of Higgs bundles and harmonic metrics (\cite{SHiggs}).

The presence of a non-trivial imaginary part $\Theta_{Im}$ seems to
imply a contradiction. Indeed, a Hodge connection
$\Theta:\Lambda^1(X,\C) \to \rho^*\Lambda^1(M,\C)$ must by
definition be compatible with the Hodge bundle structures -- in
particular, it must commute with the complex conjugation map. But
this is different from ``real''. The reason for this is the twist by
the involution $\iota:X \to X$ that we have introduced in the
definition of a Hodge bundle. This can be seen clearly if instead of
the connection $\Theta$ one considers the associated derivation
$D$. The derivation $D:\Lambda^0(X,\C) \to \rho^*\Lambda^1(M,\C)$
satisfies
\begin{equation}\label{notreal}
D\iota^*\overline{f} = \iota^*\overline{Df}
\end{equation}
for any $\C$-valued smooth function $f \in \Lambda^0(X,\C)$. Take
the decomposition $D = D_- + D_+$ into the odd and the even part
with respect to the involution $\iota$, so that every function $f
\in \Lambda^0(X,\C)$ we have
$$
D_-(\iota^*f) = -\iota^*D_-(f), \qquad\qquad 
D_+(\iota^*f) = \iota^*D_+(f),
$$
There is no reason for either one of these parts to vanish. But
\eqref{notreal} shows that
\begin{align}
D_- &= \Theta_{Re} \circ d,\label{odd}\\
D_+ &= \sqrt{-1}\Theta_{Im} \circ d,\label{even}
\end{align}
where $d:\Lambda^0(X) \to \Lambda^1(X)$ is the de Rham differential.

The imaginary part $\Theta_{Im}$ of a Hodge connection $\Theta$ on
$X/M$ -- or rather, the associated map $R_J$ -- by itself has a very
direct geometric meaning in terms of the hypercomplex structure on
$X$ given by $\Theta$. To describe it, consider the splitting
\begin{equation}\label{horvert}
\Lambda^1(X,\R) = \Lambda^1(X/M,\R) \oplus \rho^*\Lambda^1(M,\R)
\end{equation}
given by the real part $\Theta_{Re}$ and identify
$\rho^*\Lambda^1(M,\R) \cong \Lambda^1(X/M,\R)$.

\begin{lemma}\label{j.matrix}
The operator $j:\Lambda^1(X,\R) \to \Lambda^1(X,\R)$ of the
hypercomplex structure given by $\Theta$ can be written with respect
to the decomposition \eqref{horvert} as the matrix
$$
\begin{pmatrix} 0 & -R_J^{-1} \\ R_J & 0 \end{pmatrix},
$$
where $R_J:\rho^*\Lambda^1(M,\R) \to \rho^*\Lambda^1(M,\R)$ is the
bundle endomorphism \eqref{Rj.eq}.
\end{lemma}

\proof{} Since $j^2 = -\id$, it suffices to prove that for every
$1$-form $\alpha \in \rho^*\Lambda^1(M,\R)$ which lies in the
horizontal part of \eqref{horvert}, the $1$-form $j(\alpha)$ is
vertical, -- that is,
\begin{equation}\label{1.eq}
\Theta_{Re}(j(\alpha)) = 0,
\end{equation}
and moreover, that we have
\begin{equation}\label{2.eq}
\Theta_{Im}(j(\alpha)) = - \alpha.
\end{equation}
Let $\alpha$ be such a form. By definition, the kernel $\Ker\Theta
\subset \Lambda^1(X,\C)$ of the projection $\Theta:\Lambda^1(X,\C)
\to \rho^*\Lambda^1(M,\C)$ is the subbundle $\Lambda^{1,0}_J(X)$ of
$(1,0)$-forms for the complementary complex structure on
$X$. Therefore we have
$$
\alpha - \sqrt{-1}j(\alpha) \in \Ker\Theta,
$$
which means that
$$
\Theta(\alpha) = \sqrt{-1}\Theta(j(\alpha)).
$$
Since $\alpha = \Theta(\alpha)$, equations \eqref{1.eq} and
\eqref{2.eq} are the real and the imaginary parts of this
equality. \endproof

In keeping with the general philosophy of this section, we will use
the formula for $j$ given by Lemma~\ref{Rj} to express the bundle
endomorphism $R_J:\rho^*\Lambda^1(M,\C) \to \Lambda^1(M,\C)$
entirely in terms of operators on the algebra
$\rho^*\Lambda^\hdot(M,\C)$. To do this, consider the tautological
section of the pullback tangent bundle $\rho^*T(M)$, and let
$$
\tau:\rho^*\Lambda^{\hdot+1}(M,\C) \to \rho^*\Lambda^\hdot(M,\C)
$$
be the operator given by contraction with this tautological
section. Thus $\tau$ vanishes on functions, and for every $1$-form
$\alpha \in \Lambda^1(M,\R)$ the function $\tau(\rho^*\alpha)$ is
just $\alpha$ considered as a fiberwise-linear function on the total
space $\overline{T}(M)$.

\begin{lemma}\label{Rj}
For any $1$-form $\alpha \in \Lambda^1(M,\C)$ we have
$$
R_J(\alpha) = -\sqrt{-1}D_+\tau(\alpha).
$$
\end{lemma}

\proof{} By \eqref{even}, the right-hand side is equal to
$\Theta_{Im}(d\tau(\alpha)) \in \rho^*\Lambda^1(M,\C)$. Since the
projection $\Theta_{Im}:\Lambda^1(X,\C) \to \rho^*\Lambda^1(M,\C)$
vanishes on the subbundle $\rho^*\Lambda^1(M,\C) \subset
\Lambda^1(X,\C)$, this expression depends only the relative $1$-form
$P(d\tau(\alpha)) \in \Lambda^1(X/M,\C)$ obtained from the $1$-form
$d\tau(\alpha) \in \Lambda^1(X,\C)$ by the projection
$P:\Lambda^1(X,\C) \to \Lambda^1(X/M,\C)$. But $P(d\tau(\alpha))$ is
precisely the image of the form $\alpha \in \rho^*\Lambda^1(M,\C)$
under the canonical isomorphism $\rho^*\Lambda^1(M,\C) \cong
\Lambda^1(X/M,\C)$. \endproof

We will now rewrite in the same spirit the normalization condition
\eqref{norm} on the hypercomplex structure on $X$ associated to
$D$. For this we need to extends the canonical isomorphism
$\Lambda^1(X/M,\C) \cong \rho^*\Lambda^1(M,\C)$ to an algebra
isomorphism $\Lambda^\hdot(X/M,\C) \cong
\rho^*\Lambda^\hdot(M,\C)$. Then the map
$$
\tau:\rho^*\Lambda^{\hdot+1}(M,\C) \to \rho^*\Lambda^\hdot(M,\C)
$$
becomes the contraction with the relative Euler vector field (that
is, the differential of the $\R^*$-action by dilatations along the
fibers of the projection $\rho:X \to M$).

The normalization condition \eqref{norm} involves a different vector
field -- na\-me\-ly, the differential $\phi$ of the standard action
of the group $U(1)$. It will be more convenient now to multiply it
by $\sqrt{-1}$ (or, equivalently, to change the generator of the Lie
algebra of the circle $U(1)$ from $\frac{\6}{\6\theta}$ to
$z\frac{\6}{\6 z}$). Denote by
$$
\begin{CD}
\sigma:\rho^*\Lambda^{\hdot+1}(M,\C) \cong \Lambda^{\hdot+1}(X/M,\C)
@>>> \rho^*\Lambda^\hdot(M,\C) \cong \Lambda^\hdot(X/M,\C)
\end{CD}
$$
the contraction with the vertical vector field $\sqrt{-1}\phi$. The
operators $\sigma$ and $\tau$ are related by
$$
\sigma(\alpha) = \sqrt{-1}\tau(I\alpha), \qquad \alpha \in
\Lambda^1(M,\C),
$$
where $I:\Lambda^1(M,\C) \to \Lambda^1(M,\C)$ is the complex
structure operator -- in other words,
$$
\sqrt{-1}I = \begin{cases} -\id &\text{on} \quad \Lambda^{1,0}(M).\\
\id &\text{on} \quad \Lambda^{0,1}(M).\end{cases}
$$

\begin{lemma}\label{D.norm}
Let $\Theta$ be a Hodge connection on $X/M$, and let $D_+$ be the
even component of the associated derivation $D:\Lambda^0(X,\C) \to
\rho^*\Lambda^1(M,\C)$.

The hypercomplex structure on $X$ given by $\Theta$ satisfies the
normalization condition \eqref{norm} if and only for every $1$-form
$\alpha \in \Lambda^1(M,\C)$ we have
$$
\sigma \circ D_+ (f) = f.
$$
where $f = \tau(\rho^*\alpha) \in \Lambda^0(X,\C)$. Moreover, if
this is the case, then the Hodge connection $\Theta$ is holonomic.
\end{lemma}

\proof{} It suffices to check \eqref{norm} by evaluating both sides
on every $1$-form $\alpha \in \rho^*\Lambda^1(M,\C)$. Moreover, it
is even enough to check it for forms of the type $\rho^*\alpha$,
where $\alpha \in \Lambda^1(M,\C)$ is a $1$-form on $M$. Let
$\alpha$ be such a form. We have to check that
$$
j(\rho^*\alpha) \cntrct \phi = \tau(\rho^*\alpha).
$$
By Lemma~\ref{j.matrix} this is equivalent to
$$
\sigma(R_J(\rho^*\alpha)) = -\sqrt{-1}\tau(\rho^*\alpha),
$$
and by Lemma~\ref{Rj} this can be further rewritten as
$$
-\sqrt{-1}\sigma(D_+\tau(\rho^*\alpha)) =
-\sqrt{-1}\tau(\rho^*\alpha).
$$
Replacing $\tau(\rho^*\alpha)$ with $f$ gives precisely the first claim of
the lemma. 

To prove the second claim, we have to show that the map
$$
\Theta:\Lambda^1(X,\R) \to \rho^*\Lambda^1(M,\C)
$$
is surjective. By definition, on the second term
$\rho^*\Lambda^1(M,\R) \subset \Lambda^1(X,\R)$ in the splitting
\eqref{horvert} we have $\Theta = \Theta_{Re} = \id$. Therefore it
suffices to prove that
$$
\Theta_{Im}:\Lambda^1(X/M,\R) \to \sqrt{-1}\rho^*\Lambda^1(M,\R)
$$
is surjective. But by \eqref{even} and the first claim of the lemma,
this is the inverse map to 
$$
\sigma:\sqrt{-1}\rho^*\Lambda^1(M,\R) \to S^1(M,\R) \cong
\Lambda^1(X/M,\R). \qquad\qquad\qquad\square
$$

\section{The Weil algebra.}\label{wa}

The final preliminary step in the proof of Theorem~\ref{hc} is to
reduce it from a question about the total space $X = \overline{T}M$
of the complex-conjugate to the tangent bundle on $M$ to a question
about the manifold $M$. To do this. we introduce the following.

\begin{defn}\label{wa.def}
The {\em Weil algebra} $\B^\hdot(M)$ of a complex manifold $M$ is
the algebra on $M$ defined by
$$
\B^\hdot(M) = \rho_*\rho^*\Lambda^\hdot(M,\C),
$$
where $\rho:\overline{T}M \to M$ is the canonical projection.
\end{defn}

This requires an explanation -- indeed, for a vector bundle $\E$ on
$\overline{T}M$, the direct image sheaf $\rho_*\E$ {\em a priori} is
not a sheaf of sections of any vector bundle on $M$. We have to
consider a smaller subsheaf. From now on and until Theorem~\ref{hc}
is proved, we will be interested not in hypercomplex structures on
the total space $\overline{T}M$ but in their formal Taylor
decompositions in the neighborhood of the zero section $M \subset
\overline{T}M$. Therefore it will be sufficient for our purposes to
define the direct image $\rho_*\E$ as the sheaf of sections of the
bundle $\E$ on $\overline{T}M$ {\em which are polynomial along the
fibers of the projection $\rho:\overline{T}M \to M$}. Formal germs
of bundle maps on $\E$ will give formal series of maps between the
corresponding direct image bundles.

Having said this, we can explicitly describe the Weil algebra
$\B^\hdot(M)$. Our first remark is that $\B^k(M)$ is canonically a
Hodge bundle on $M$ of weight $k$. Moreover, since the $U(1)$-action
on $M$ is trivial, Hodge bundles on $M$ are just bundles of
$\R$-Hodge structures in the usual sense. Thus we have a Hodge type
bigrading
$$
\B^k(M) = \bigoplus_{p+q=k} \B^{p,q}(M)
$$
and a canonical real structure on every one of the complex vector
bundles $\B^k(M)$.

The projection formula show that for every $k$ we have a canonical
isomorphism 
$$
\B^k(M) \cong \B^0(M) \otimes \Lambda^k(M,\C).
$$
These isomorphisms are compatible with the Hodge structures and with
multiplication. The degree-$0$ Hodge bundle $\B^0(M)$ is a symmetric
algebra freely generated by the bundle $S^1(M,\C)$ of functions on
$\overline{T}M$ linear along the fibers of $\rho:\overline{T}M \to
M$. The complex vector bundle $S^1(M,\C)$ is canonically isomorphic
to the bundle $\Lambda^1(M,\C)$ of $1$-forms on $M$. However, the
Hodge structures on these bundles are different. The Hodge type
grading on $S^1(M,\C)$ is given by
$$
S^1(M,\C) = S^{1,-1}(M) \oplus S^{-1,1}(M,\C),
$$
where $S^{1,-1}(M) \cong \Lambda^{1,0}(M)$ and $S^{-1,1}(M) \cong
\Lambda^{0,1}(M)$ -- the grading is the same as on $\Lambda^1(M,\C)$
but graded pieces are assigned different weights. Moreover, the
complex conjugation map on $S^1(M,\C)$ is {\em minus} the complex
conjugation map on $\Lambda^1(M,\C)$. This is the last vestige of
the twist by the involution $\iota:\overline{T}M \to \overline{T}M$
in Definition~\ref{hb.def}.

To simplify notation, denote by $S^k(M,\C)$ the $k$-th symmetric
power of the Hodge bundle $S^1(M,\C)$. Then we have
$$
\B^0(M) = \bigoplus_{k \geq 0}S^k(M,\C),
$$
and the Weil algebra $\B^\hdot(M) = \B^0(M) \otimes
\Lambda^\hdot(M,\C)$ is the free graded-commutative algebra
generated by $S^1(M,\C)$ and $\Lambda^1(M,\C)$ (where, contrary to
notation, $S^1(M,\C)$ is placed in degree $0$).

It will be convenient to introduce another grading on the Weil
algebra $\B^\hdot(M)$ by assigning to both of the generator bundles
$S^1(M,\C)$, $\Lambda^1(M,\C)$ degree $1$. We will call it {\em
augmentation grading} and denote by lower indices, so that we have
\begin{align*}
S^1(M,\C) &= \B^0_1(M) \subset \B^0(M),\\ 
\Lambda^1(M,\C) &= \B^1_1(M) \subset \B^1(M).
\end{align*}
The augmentation grading corresponds to the Taylor decomposition
near the zero section $M \subset \overline{T}M$. Namely, every
formal germ near $M \subset \overline{T}M$ of a flat Hodge
connection on $\overline{T}M/M$ induces a formal series
\begin{equation}\label{taylor}
D = \sum_{k \geq 0} D_k
\end{equation}
of algebra bundle derivations
$$
D_k:\B^\hdot(M) \to \B^{\hdot+1}(M),
$$
where each of the derivations $D_k$ is weakly Hodge and has
augmentation degree $k$. Their (formal) sum satisfies
$$
D \circ D = 0.
$$
Conversely, every formal series \eqref{taylor} induces a (formal
germ of a) weakly Hodge derivation $D:\rho^*\Lambda^\hdot(M,\C) \to
\rho^*\Lambda^{\hdot+1}(M,\C)$ on the total space
$\overline{T}M$. This derivation comes from a flat Hodge connection
if and only if we have
$$
D(f) = df
$$
for every function $f \in \B^0_0(M) \cong \Lambda^0(M,\C)$. Since we
have $D \circ D = 0$, this immediately implies that $D$ coincides
with the de Rham differential $d$ on the whole subalgebra
$$
\Lambda^\hdot(M,\C) \subset \B^\hdot(M) \cong \B^0(M) \otimes
\Lambda^\hdot(M,\C).
$$
More precisely, we must have $D_1 = d$ on $\Lambda^\hdot(M,\C)$, and
all the other components $D_k, k \neq 1$ must vanish on this
subalgebra. Since $D$ is a derivation, this in turn implies that all
the components $D_k:\B^\hdot(M) \to \B^{\hdot + 1}(M)$ for $k \neq
1$ are not differential operators but bundle maps. Moreover, since
the algebra $\B^\hdot(M)$ is freely generated by $S^1(M,\C)$ and
$\Lambda^1(M,\C)$, and we know {\em a priori} the derivation $D$ on
the generator subbundle $\Lambda^1(M,\C)$, it always suffices to
specify the restriction $D:S^1(M,\C) \subset \B^0(M) \to \B^1(M)$.

The decomposition $D = D_- + D_+$ of a Hodge connection into an even
and an odd part is quite transparent on the level of the Weil
algebra -- we simply have
$$
D_- = \sum_{k \geq 0}D_{2k+1} \qquad\qquad 
D_+ = \sum_{k \geq 0}D_{2k}.
$$
We will now rewrite the normalization condition \eqref{norm} in
terms of the Weil algebra. To do this, note that the map
$\sigma:\rho^*\Lambda^{\hdot+1}(M,\C) \to \Lambda^\hdot(M,\C)$
induces a bundle map $\sigma:\B^{\hdot+1}(M) \to \B^\hdot(M)$. This
map is in fact a derivation of the Weil algebra. It vanishes on the
generator bundle $S^1(M,\C)$, while on the generator bundle
$\Lambda^1(M,\C)$ it is given by
$$
\sigma = \begin{cases} \id:\Lambda^{1,0}(M) \to S^{1,-1}(M) \cong
\Lambda^{1,0}(M), \\ -\id:\Lambda^{0,1}(M) \to S^{-1,1}(M) \cong
\Lambda^{0,1}(M). \end{cases}
$$
Then Lemma~\ref{D.norm} immediately shows that the (formal germ of
the) hypercomplex structure on $\overline{T}M$ induced by a
derivation $D:\B^\hdot(M) \to \B^{\hdot+1}(M)$ is normalized if and
only if we have
\begin{equation}\label{norm.bis}
\sigma \circ D_+ = \id
\end{equation}
on the generator subbundle $S^1(M,\C) \subset \B^0(M)$. It is
convenient to modify this in the following way. Let $C:S^1(M,\C) \to
\Lambda^1(M,\C)$ be the isomorphism inverse to
$\sigma:\Lambda^1(M,\C) \to \Lambda^1(M,\C)$. Set $C = 0$ on the
generator subbundle $\Lambda^1(M,\C) \subset \B^1(M)$ and extend it
a derivation $C:\B^\hdot(M) \to \B^{\hdot+1}(M)$ of the Weil
algebra. Both derivations $C$ and $\sigma$ are real. Moreover, the
derivation $C$ is weakly Hodge (the derivation $\sigma$ is not --
simply because it decreases the weight). Then the normalization
condition is equivalent to
$$
\begin{cases} D_0 = C,\\ \sigma \circ D_k = 0 \text{ on } S^1(M,\C)
\subset \B^0(M) \text{ for every even } k = 2p \geq 1.
\end{cases}
$$
To sum up, formal germs near $M \subset \overline{T}M$ of normalized
flat Hodge connections on $\overline{T}M$ are in a natural
one-to-one correspondence with derivations
$$
D = \sum_{k \geq 0} D_k:\B^\hdot(M) \to \B^{\hdot+1}(M)
$$
of the Weil algebra $\B^\hdot(M)$ which satisfy the following
conditions.
\begin{enumerate}
\item $D \circ D = 0$.
\item $D_k:\B^\hdot(M) \to \B^{\hdot+1}(M)$ is a weakly Hodge
algebra derivation of augmentation degree $k$.
\item $D_0 = C$.
\item $D_1 = d$ and $D_k = 0$, $k \neq 0$ on the subalgebra
$\Lambda^\hdot(M,\C) \subset \B^\hdot(M)$.
\item $\sigma \circ D_{2k} = 0$ on $S^1(M,\C) = \B^0_1(M) \subset
\B^0(M)$ for every $k \geq 1$.
\end{enumerate}\label{cnd}
For every such derivation, the differential operator 
$$
D_1:S^1(M,\C) = \B^0_1(M) \to \B^1_2(M) \cong S^1(M,\C) \otimes
\Lambda^1(M,\C) 
$$
satisfies the Leibnitz rule
$$
D_1(fa) = fD_1(a) + adf, \qquad a \in S^1(M,\C), f \in
\Lambda^0(M,\C).
$$
Therefore it is a connection on the bundle $S^1(M,\C)$. We postpone
the proof of the following Lemma till the end of
Section~\ref{metrics}.

\begin{lemma}\label{obata}
The connection $D_1$ on the bundle $S^1(M,\C) \cong \Lambda^1(M,\C)$
coincides with the connection on $M$ induced by the Obata connection
for the hypercomplex structure on $\overline{T}M$ defined by the
derivation $D$.
\end{lemma}

With Lemma~\ref{obata} in mind, we see that Theorem~\ref{hc} is
reduced to the following statement.

\begin{prop}\label{hc.prop}
Let $\nabla$ be a torsion-free connection on the cotangent bundle
$\Lambda^1(M,\C)$ of a complex manifold $M$. Assume that the
curvature of the connection $\nabla$ is of type $(1,1)$.

Then there exists a unique derivation 
$$
D = \sum_{k \geq 0}D_k:\B^\hdot(M) \to \B^{\hdot+1}(M)
$$ 
of the Weil algebra $\B^\hdot(M)$ of the manifold $M$ such that $D$
satisfies the conditions \thetag{i}-\thetag{v} above and we have
$$
D_1 = \nabla
$$
on $S^1(M,\C) \subset \B^0(M)$.
\end{prop}

This ends the preliminaries. We now begin the proof of
Proposition~\ref{hc.prop}.

\section{The proof of Proposition~\latexref{hc.prop}.}\label{pf}

The proof proceeds by induction on the augmentation degree. Denote
$$
D_{\leq k} = D_0 + D_1 + \dots + D_k.
$$
To base the induction, consider the derivation $D_{\leq 1}$. By
assumptions it is equal to
$$
D_{\leq 1} = C + D_1.
$$
Let $R:\B^\hdot(M) \to \B^{\hdot+1}(M)$ be the composition
$$
R = D_{\leq 1} \circ D_{\leq 1}.
$$
Since the derivation $D_{\leq 1}$ of the graded-commutative algebra
$\B^\hdot(M)$ is of odd degree, up to a coefficient the composition
$R$ coincides with the supercommutator $\{D_{\leq 1},D_{\leq
1}\}$. In particular, it is also an algebra derivation.

The derivation $R:\B^\hdot(M) \to \B^{\hdot+2}(M)$ {\em a priori}
has components $R_0$, $R_1$, $R_2$ of augmentation degrees $0$, $1$
and $2$. However,
$$
R_0 = \{C,C\} = 0.
$$
Moreover,
$$
R_1 = \{C,D_1\}:S^1(M,\C) \to \B^2_2(M) \cong \Lambda^2(M,\C)
$$
is precisely the torsion of the connection $\nabla = D_1$. Thus it
vanishes by assumption. Since $C = 0$ on $\Lambda^\hdot(M,\C)
\subset \B^\hdot(M)$, we also have $R_1 = 0$ on $\Lambda^1(M,\C)$,
which implies that $R_1 = 0$ everywhere. What remains is $R_2$. In
general, it does {\em not} vanish, and to kill it we have to add new
terms $D_k$.

We now turn to the induction step. Assume that for some $k \geq 2$
we are already given the derivation $D_{\leq k-1}$ which satisfies
the conditions \thetag{ii}-\thetag{v} on page \pageref{cnd}, and
assume that the composition $D_{\leq k-1} \circ D_{\leq k-1}$ has no
non-trivial components of augmentation degrees $< k$. Denote by
$R_k:\B^\hdot(M) \to \B^{\hdot+1}(M)$ its component of augmentation
degree $k$. We have to find a derivation $D_k:\B^\hdot(M) \to
\B^{\hdot+1}$ which also satisfies \thetag{ii}-\thetag{v} and such
that
$$
(D_{\leq k-1} + D_k) \circ (D_{\leq k-1} + D_k) = 0
$$
in augmentation degree $k$. This is equivalent to
\begin{equation}\label{Rk}
\{D_0,D_k\} = \{C,D_k\} = -R_k.
\end{equation}
The conditions \thetag{ii} and \thetag{iii} mean that
$D_k:\B^\hdot(M) \to \B^{\hdot+1}(M)$ must be a weakly Hodge
derivation which vanishes on $\Lambda^1(M,\C) = \B^1_1(M) \subset
\B^1(M)$. The condition \thetag{iv} is only relevant for
$D_1$. Finally, the condition \thetag{v} is relevant for all even
$k$ and means that
$$
\sigma \circ D_k = 0
$$
on $S^1(M,\C) \subset \B^0(M)$.

Because of \thetag{iii}, it suffices to define $D_k$ on the
generator subbundle $S^1(M,\C)$. Since $D_k$ commutes with the
complex conjugation map, it even suffices to consider only
$S^{1,-1}(M) \subset S^1(M,\C)$. Moreover, \thetag{iii} implies that
it also suffices to check \eqref{Rk} only on $S^{1,-1}(M) =
\B^{1,-1}_1(M) \subset \B^0(M)$. Note that on this subbundle we have
$\{C,D_k\} = C \circ D_k$.

There will be two slightly different cases. The first is one when
$k=2p+1$ is odd, the second one is when $k = 2p$ is even.

In both cases, the weakly Hodge map $R_k:\B^0_1(M) \to
\B^2_{k+1}(M)$ only has non-trivial pieces of Hodge bidegrees
$(2,0)$, $(1,1)$ and $(0,2)$. Moreover, for any map $\Theta$ of odd
degree we tautologically have $[\{\Theta,\Theta\},\Theta] =
0$. Applying this to $\Theta = D_{\leq k}$ and collecting terms of
augmentation degree $k$, we see that
$$
C \circ R_k = 0.
$$
To track various components of the map $R_k:\B^0_1(M) \to
\B^2_{k+1}(M)$ it is convenient to refine the augmentation grading
on the Weil algebra $\B^\hdot(M)$ to an {\em augmentation bigrading}
by setting
\begin{align*}
\deg S^{1,-1}(M) &= \deg \Lambda^{1,0}(M) = (1,0),\\
\deg S^{-1,1}(M) &= \deg \Lambda^{0,1}(M) = (0,1),
\end{align*}
on the generator bundles $S^1(M,\C) = S^{1,-1}(M) \oplus
S^{-1,1}(M)$ and $\Lambda^1(M,\C) = \Lambda^{1,0}(M) \oplus
\Lambda^{0,1}(M)$. The augmentation bigrading will be denoted by
lower indices, so that we have
$$
\B^\hdot_k = \bigoplus_{p+q = k}\B^\hdot_{p,q}.
$$
The relevant pieces of the augmentation bigrading on the bundle
$\B^\hdot_{k+1}(M)$ are shown on Figure~\ref{fig.odd} for $k = 2p+1$
odd, and on Figure~\ref{fig.even} for $k = 2p$ even. The axes on the
figures correspond to the grading by Hodge type. A Hodge bidegree
component $\B^{p,q}_{m,n}$ can be non-trivial only when $p \geq m-n$
and $q \geq n-m$. Thus the component $\B^{\hdot,\hdot}_{m,n}$ is
represented by an upward-looking angle with vertex $(m-n,n-m)$: a
graded piece $\B^{p,q}_{m,n}$ can be non-trivial only if the point
$(p,q)$ lies in the interior (or on the boundary) of this angle.

\begin{figure}[t]
\epsfig{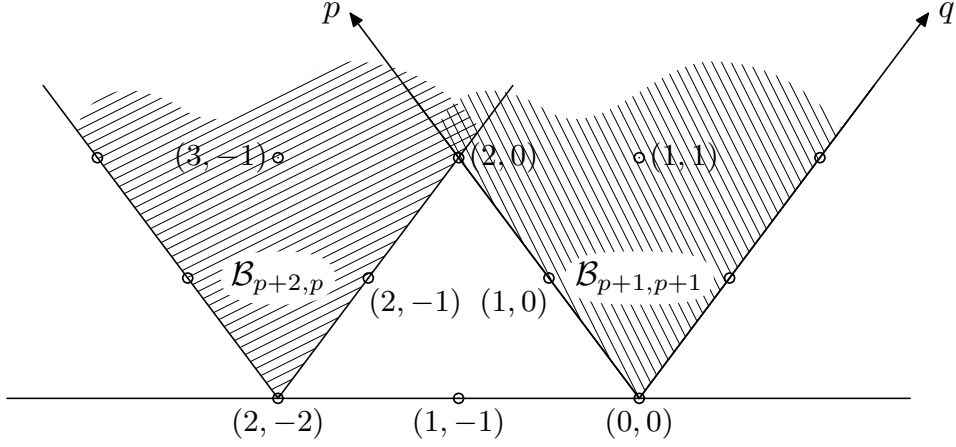}
\caption{The augmentation bigrading on
$\protect\B^\protect\hdot_{k+1}$ for an odd $k$, $k = 2p+1$.}
\end{figure}\refstepcounter{figuredammit}\label{fig.odd}

\begin{figure}[t]
\epsfig{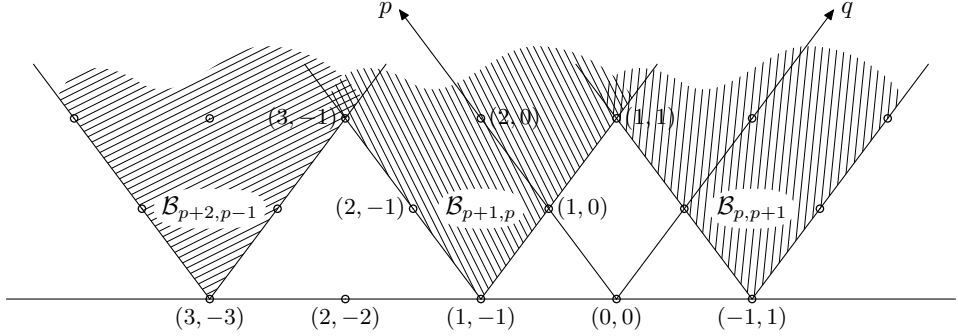}
\caption{The augmentation bigrading on
$\protect\B^\protect\hdot_{k+1}$ for an even $k$, $k = 2p$.}
\end{figure}\refstepcounter{figuredammit}\label{fig.even}

Consider the Hodge bidegree decompositions
$$
C = C^{1,0} + C^{0,1} \qquad\qquad \sigma = \sigma^{-1,0} +
\sigma^{0,-1}
$$
of the derivations $C$, $\sigma$ of the Weil algebra
$\B^\hdot(M)$. Then the augmentation bigrading is essentially the
eigenvalue decomposition for the commutators
$\{C^{1,0},\sigma^{-1,0}\}$ and $\{C^{0,1},\sigma^{0,-1}\}$. More
precisely, we have
\begin{align}\label{com.aug}
\begin{split}
\{C^{1,0},\sigma^{0,-1}\} &= \{C^{0,1},\sigma^{-1,0}\} = 0,\\
\{C^{1,0},\sigma^{-1,0}\} &= m\id \text{ on }\B^\hdot_{m,n},\\
\{C^{0,1},\sigma^{0,-1}\} &= n\id \text{ on }\B^\hdot_{m,n}
\end{split}
\end{align}
Indeed, since all these commutators are derivations of the Weil
algebra, it suffices to check this on the generator bundles
$S^1(M,\C)$, $\Lambda^1(M,\C)$, which is elementary. In particular,
we see that both $C$ and $\sigma$ preserve the augmentation
bidegree. The equalities \eqref{com.aug} also immediately imply that
$$
\{C,\sigma\} = k\id \text{ on } \B^\hdot_k.
$$
One further corollary of \eqref{com.aug} will be very important (we
leave the proof to the reader as an easy exercise). 

\begin{lemma}\label{bndry}
If $m,n \geq 1$, then the map $C^{1,0}$ is injective on every graded
piece $\B^{p,q}_{m,n}$ with $q = n-m$, while $C^{0,1}$ is injective
on $\B^{p,q}_{m,n}$ with $p=m-n$.
\end{lemma}

Graphically, this means that $C^{1,0}$ is injective on
$\B^{p,q}_{m,n}$ when the point $(p,q)$ lies on the right-hand
boundary of the angle representing $\B^{\hdot,\hdot}_{m,n}$, and
$C^{0,1}$ is injective in this graded piece when the point $(p,q)$
lies on the left-hand boundary of the same angle. We will call this
{\em the boundary rule}.

We can now proceed with the proof of the induction step.

\medskip

\noindent{\em Case when $k=2p+1$ is odd.} 
Looking at Figure~\ref{fig.odd}, we see that the only non-trivial
augmentation-bidegree components of the map $R_k:\B^0_1(M) \to
\B^2_{k+1}(M)$ are $R_{p,p+1}$ and $R_{p+1,p}$,
$$
R_k = R_{p,p+1} + R_{p+1,p},
$$
and the same is true for any weakly Hodge map $D_k:S^{1,-1}(M) \to
\B^1_{k+1}(M)$,
$$
D_k = D_{p,p+1} + D_{p+1,p}.
$$
Moreover, on $S^{1,-1}(M) \subset \B^0(M)$ we have
\begin{align*}
D^{0,1}_k &= D_{p,p+1} \qquad\qquad D^{1,0}_k = D_{p+1,p},\\
R^{0,2}_k &= R^{0,2}_{p,p+1} \qquad\qquad D^{2,0}_k = R^{2,0}_{p+1,p},
\end{align*}
while $R_k^{1,1}$ further decomposes as $R^{1,1}_{p+1,p} +
R^{1,1}_{p,p+1}$. 

We have to find $D_k$ which satisfies \eqref{Rk}. In particular, we
must have
\begin{equation}\label{first.line}
R^{2,0}_k = -C^{1,0} \circ D^{1,0}_k \qquad\qquad R^{0,2}_k = -
C^{0,1} \circ D^{0,1}_k.
\end{equation}
But by the boundary rule the map $C^{1,0}$ is injective on
$\B^{2,-1}_{p+2,p}$, while the map $C^{0,1}$ is injective on
$\B^{1,0}_{p+1,p+1}$. Therefore there exists at most one weakly
Hodge map $D_k:\B^0_{1,0}(M) \to \B^1_{k+1}$ which satisfies
\eqref{first.line}. Setting (again on $S^{1,-1}(M)$)
\begin{align*}
D^{0,1}_k &= D^{0,1}_{p,p+1} = -\frac{1}{p+1}\sigma^{0,-1} \circ
R^{0,2}_k,\\ 
D^{1,0}_k &= D^{1,0}_{p+1,p} = -\frac{1}{p+2}\sigma^{-1,0} \circ
R^{2,0}_k. 
\end{align*}
gives this unique solution to \eqref{first.line}. Indeed, we have
\begin{align*}
C^{1,0} \circ D^{1,0}_k &= -\frac{1}{p+2} C^{1,0} \circ \sigma
^{-1,0} \circ R^{2,0}_k \\ 
&= -\frac{1}{p+2}\left(\{C^{1,0},\sigma^{-1,0}\}\circ
R^{2,0}_{p+1,p} + \sigma^{-1,0} \circ C^{1,0} \circ
R^{2,0}_k\right).
\end{align*}
The second summand in the brackets vanishes since $C \circ R_k = 0$,
while the first is equal to $(p+2)R^{2,0}_k$ by
\eqref{com.aug}. This proves the first equation in
\eqref{first.line}. The second one is proved in exactly the same
way.

It remains to prove that this map $D_k$ satisfies not only
\eqref{first.line} but also the stronger condition \eqref{Rk}. To do
this, note that
$$
C \circ (C \circ D_k + R_k) = (C \circ C) \circ D_k + C \circ R_k =
0.
$$
But from \eqref{first.line} we see that $C \circ D_k + R_k$ is of
Hodge bidegree $(1,1)$. Therefore this implies that
$$
C^{1,0} \circ (C \circ D_k + R_k) = C^{0,1} \circ (C \circ D_k +
R_k) = 0.
$$
The only possible non-trivial components of $C \circ D_k + R_k$ with
respect to the augmentation bigrading have bidegrees $(p+1,p)$ and
$(p,p+1)$, and by the boundary rule $C^{1,0}$ is injective on
$\B^{2,0}_{p+1,p}$, while $C^{0,1}$ is injective on
$\B^{2,0}_{p,p+1}$. Thus $C \circ D_k + R_k = 0$.

\medskip

\noindent{\em Case when $k=2p$ is odd.} 
Assume for the moment that $k \geq 4$, thus $p \geq 2$.

Looking at Figure~\ref{fig.even}, we see that {\em a priori} the map
$R_k$ can have three non-trivial augmentation-bidegree components,
namely, 
$$
R_k = R^{0,2}_{p-1,p+1} + R_{p,p} + R^{2,0}_{p+1,p-1}.
$$
However, since $C \circ R_k = 0$, and the map $C^{1,0}$ is injective
on $\B^{3,-1}_{p+2,p-1}(M)$ by the boundary rule, we see that
$R^{0,2}_{p+1,p-1} = 0$. Analogously, $R_{p-1,p+1} = 0$. Therefore
$R_k$ is of pure augmentation bidegree $(p,p)$.

Since the map $D_k:S^1(M) \to \B^1_{k+1}(M)$ is weakly Hodge, it
must also be of augmentation bidegree $(p,p)$. Conversely, looking
at the angle representing $\B^{\hdot,\hdot}_{p+1,p-1}$, we see that
{\em every} real map $D_k:S^1(M,\C) \to \B^1_{k+1}$ of pure
augmentation bidegree $(p,p)$ is necessarily weakly Hodge. In
particular, setting
\begin{equation}\label{even.k}
D_k = -\frac{1}{k+1}\sigma \circ R_k
\end{equation}
defines a weakly Hodge map. This map is a solution to \eqref{Rk}:
$$
C \circ D_k = -\frac{1}{k+1}C \circ \sigma \circ R_k =
-\frac{1}{k+1}\{C,\sigma\} \circ R_k - \sigma \circ C \circ R_k =
-R_k.
$$
This solution is not unique. However, since $k$ is even, we have the
additional normalization condition $\sigma \circ D_k = 0$. This
condition (automatically satisfied by the solution \eqref{even.k})
ensures uniqueness. Indeed, the difference $P = D_k - D_k'$ between
two solutions $D_k$, $D_k'$ must satisfy $C \circ P = \sigma \circ P
= 0$, which implies
$$
P = \frac{1}{k+1}\{C,\sigma\} \circ P = 0.
$$
Finally, it remains to consider the case $k=2$. The general argument
works in this case just as well, with a single exception. Since $p-1
= 0$ is no longer strictly positive, the boundary rule does not
apply: it is not true that $C^{1,0}$ is injective on
$\B^{3,-1}_{p+1,p-1} = \B^{3,-1}_{3,0}$ (in fact, on this graded
piece $C^{1,0}$ is equal to zero). Therefore the component
$R_{p+1,p-1}$ does not vanish automatically. However, this component
$$
R^{2,0}_{2,0}:S^{1,-1}(M) \to \B^{3,-1}_{3,0}(M) \cong
\Lambda^{2,0}(M) \otimes S^{1,-1}(M)
$$
is precisely the $(2,0)$-curvature of the connection $\nabla$ on
$M$. It vanishes by the second assumption on this
connection. \endproof

\section{Metrics.}\label{metrics}

The last Section essentially finishes the proof of the hypercomplex
Theorem~\ref{hc} (it remains to prove Lemma~\ref{obata}). We will
now sketch a proof of the hyperk\"ahler Theorem~\ref{hk}.

As we have already noted, Theorem~\ref{hk} will be a corollary of
Theorem~\ref{hc}. Namely, given a K\"ahler manifold $M$ we proceed
in the following way. First we note that the Levi-Civita connection
$\nabla_{LC}$ on $M$ has no torsion and no
$(2,0)$-curvature. Therefore Theorem~\ref{hc} applies to
$\nabla_{LC}$ and provides a hypercomplex structure on the total
space $X = \overline{T}M$. Then we show that every Hermitian metric
on $M$ which is preserved by $\nabla_{LC}$ (in particular, the given
K\"ahler metric) extends uniquely from the zero section $M \subset
\overline{T}M$ to a (formal germ of a) hyper-hermitian metric on the
hypercomplex manifold $X = \overline{T}M$ which is compatible with
the hypercomplex structure. After this, we finish the proof by
identifying the holomorphically symplectic manifolds $\overline{T}M$
and $T^*M$.

We will go through these steps in reverse order, starting with the
last one.

\begin{lemma}\label{cotang}
Assume given a hypercomplex structure on the total space $X =
\overline{T}M$ which satisfies the conditions of
Theorem~\ref{hc}. Let $h$ be a $U(1)$-invariant hyperk\"ahler metric
on $X$ compatible with this hypercomplex structure, and let $\Omega_X
\in \Lambda^{2,0}(X)$ be the associated holomorphic $2$-form. Let
$T^*M$ be the total space of the cotangent bundle to $M$ equipped
with the standard holomorphic $2$-form $\Omega$.

Then there exists a unique $U(1)$-equivariant biholomorphic map
$\eta:X \to T^*M$ such that $\Omega_X = \eta^*\Omega$.
\end{lemma}

\proof{} Since the map $\eta$ must be $U(1)$-equivariant, it must
commute with the canonical projections $\rho:X,T^*M \to M$ and send
the zero section $M \subset X$ to the zero section $M \subset
T^*M$. Denote by $\phi$ the differential of the $U(1)$-action. Then
we also must have
$$
\Omega_X \cntrct \phi = \eta^*\Omega \cntrct \phi = \eta^*(\Omega
\cntrct \phi).
$$
But the $1$-form $\alpha = \Omega \cntrct \phi$ is the
tautological $1$-form $\alpha \in \rho^*(\Lambda^1(M)) \subset
\Lambda^1(T^*M)$. Therefore the $1$-form $\rho^*\alpha$ on $X$
completely defines the map $\eta$. 

Conversely, the form $\alpha_X = \Omega_X \cntrct \phi$ satisfies
$\alpha_X = \eta^*\alpha$ for a unique map $\eta:X \to T^*M$. Since
the metric $h$ is $U(1)$-invariant, the forms $\Omega_X$ and
$\alpha_X$ are of weight $1$. Therefore the map $\eta:X \to T^*M$ is
$U(1)$-equivariant. By the Cartan homotopy formula, we have
$$
\Omega_X = d\alpha_X = d\eta^*\alpha = \eta^*d\alpha =
d\Omega. \qquad \square
$$

We will now explain how to construct the metric $h$ -- or,
equivalently, the associated holomorphic $2$-form $\Omega_J \in
\Lambda^{2,0}(X)$.

Keep the notation of last two Sections. Let $\omega \in
\Lambda^{1,1}(M,\C)$ be the K\"ahler form (more generally, any
$(1,1)$-form preserved by the connection). We have to prove that
there exists a unique (formal germ of a) holomorphic $(2,0)$-form
$\Omega \in \Lambda^{2,0}_J(X)$ which is of $H$-type $(1,1)$ and
whose restriction to the zero section $M \subset X$ coincides with
$\omega$ (since the positivity \eqref{posit} is an open condition,
it is satisfied automatically in a neighborhood of the zero section
$M \subset X$).

To reformulate this in terms of $M$, consider the complex
$\Lambda^{2,\hdot}_J(X)$ of Hodge bundles on $\overline{T}M$ with
the Dolbeault differential $D=\bar\6_J:\Lambda^{2,\hdot}_J(X) \to
\Lambda^{2,\hdot+1}(X)$. Denote by
$$
\CC^{\hdot+2}(M) = \rho_*\Lambda^{2,\hdot}_J(X)
$$
its direct image on $M$ (the grading is shifted by $2$ to make it
compatible with the Hodge degrees). We are given a section $\omega
\in \Lambda^{1,1}(M)$ of Hodge type $(1,1)$. We have to prove that
there exists a section $\Omega = \Omega_J \in \CC^{1,1}(M)$ such
that $\Omega = \omega$ on the zero section $M \subset X$ and
$D\Omega = 0$.

Since $\Lambda^{2,\hdot}_J(X) \cong \Lambda^{2,0}_J(X) \otimes
\Lambda^{0,\hdot}_J(X)$ and $\Lambda^{2,0}_J(X) \cong
\rho^*\Lambda^2(M,\C)$, the complex $\CC^\hdot(M)$ is a free module
$$
\CC^\hdot(M) \cong L^2(M) \otimes \B^\hdot(M)
$$
over the Weil algebra $\B^\hdot(M) = \rho_*\Lambda^{0,\hdot}_J(X)$
generated by some subbundle 
$$
L^2(M) \subset \CC^2(M)
$$ 
which is isomorphic to $\Lambda^2(M,\C)$. We introduce the
augmentation grading on the $\B^\hdot(M)$-module $\CC^\hdot(M)$ by
setting $\deg L^2(M) = 2$. Just as in Proposition~\ref{hc.prop}, the
proof will proceed by induction on the augmentation degree --
namely, we will construct the form $\Omega \in \CC^{1,1}(M)$ as a
sum
$$
\Omega = \Omega_0 + \dots + \Omega_k + \dots
$$
with $\Omega_k \in \CC^{1,1}_k(M)$ of augmentation degree $k+2$. We
begin with the induction step. It is completely parallel to
Proposition~\ref{hc.prop}, so we give only a sketch.

\proof[Induction step -- a sketch.] We can assume that we already
have 
$$
\Omega_{<k} = \Omega_0 +\dots+\Omega_{k-1}
$$ 
such that $\Phi = D\Omega_{<k}$ is of augmentation degree $\geq
k+2$. Denote by $\Phi_k = \Phi_k^{2,1} + \Phi_k^{1,2}$ the component
of augmentation degree exactly $k+2$. We have to show that there
exists a unique $\Omega_k \in \CC^{1,1}_{k+2}$ such that $\Phi_k =
D_0\Omega_k$.

The derivations $C$ and $\sigma$ of the Weil algebra $\B^\hdot(M)$
extend to endomorphisms of the free module $L^2 \otimes \B^\hdot(M)$
by setting $C = \sigma = 0$ on $L^2(M)$. Just as on $\B^\hdot(M)$,
we have $D_0 = C$. For every $k \geq 0$, we have $\{C,\sigma\} =
k\id$ on $\CC^\hdot_{k+2}$. This immediately implies that $C = D_0$
is injective on $\CC^{1,1}_{k+2}(M)$ for $k \geq 1$, which proves
the uniqueness of $\Omega_k$.

The space $\CC^\hdot_{k+2}$ splits into the sum of parts of the form
$$
L^{p,q} \otimes \B^\hdot_{m,n}, \qquad\qquad p+q = 2;\;m+n=k;\;p,q,m,n
\geq 0.
$$
Such a part can have a non-trivial piece of Hodge bidegree $p_1,q_1$
only if $p_1 \geq p+m-n$ and $q_1 \geq q+n-m$. Having in mind the
graphical representation as in Figure~\ref{fig.odd} and
Figure~\ref{fig.even}, we will say that the part $L^{p,q} \otimes
\B^\hdot_{m,n} \subset \CC^\hdot_{k+2}$ is {\em an angle based at
$(p+m-n,q+n-m)$}. Each angle is preserved by the maps $C$ and
$\sigma$. In this terminology, $\CC^{2,1}_{k+2}$ can intersect
non-trivially with various angles based at $(2,0)$ and $(1,1)$,
while $\CC^{1,2}_{k+2}$ can intersect non-trivially with angles
based at $(1,1)$ and $(0,2)$. But by induction we have $C\Phi_k =
0$, which implies that $C^{1,0}\Phi_k^{2,1} = C^{0,1}\Phi_k^{1,2} =
0$. Applying the boundary rule Lemma~\ref{bndry}, we see that both
$\Phi_k^{2,1} \in \CC^{2,1}_{k+2}$ and $\Phi_k^{1,2} \in
\CC^{1,2}_{k+2}$ must lie entirely within angles based at
$(1,1)$. Therefore $\sigma\Phi_k$ must also lie within angles based
in $(1,1)$. This means that $\sigma\Phi_k \in \C^{1,1}_{k+2}$ is of
Hodge type $(1,1)$, and we can set $\Omega_k =
\frac{1}{k}\sigma\Phi_k$. \endproof

We now have to base the induction -- namely, to find the section
$\Omega_0 \in \CC^{1,1}_2(M)$ with correct restriction to the zero
section $M \subset X$, and to handle those angles $L^{p,q}(M)
\otimes \B^\hdot_{m,n}$ to which the boundary rule does not apply --
which means the angles with $m=0$ or $n=0$. There are three such
angles. We denote the corresponding components of the derivation
$D:L^{1,1}(M) \to L^2(M) \otimes \B^1(M)$ by
\begin{align}
D_1:L^{1,1}(M) &\to L^{1,1}(M) \otimes \B^1_1(M),\label{d1}\\
D_2^l:L^{1,1}(M) &\to L^{2,0}(M) \otimes
\B^{-1,2}_{0,2}(M),\label{d2r}\\ 
D_2^l:L^{1,1}(M) &\to L^{0,2}(M) \otimes
\B^{2,-1}_{2,0}(M).\label{d2l} 
\end{align}
We have to choose $\Omega_0$ so that $D_1\Omega_0 = D_2^l\Omega_0 =
D_2^r\Omega_0 = 0$. We note that $D_2^l\Omega_0 = 0$ implies
$D_2^r\Omega_0 = 0$ by complex conjugation.

Moreover, we note that we have one more degree of freedom. So far
nothing depended on the choice of the generator subbundle $L^2(M)
\subset \CC^2(M)$ -- all we needed was to know that it exists. We
will now make this choice. It will not be the most obvious one, but
the one which will make computations as easy as possible. We
consider the splitting
\begin{equation}\label{biTheta}
\Lambda^1(X,\C) = \Lambda^{1,0}_J(X) \oplus \rho^*\Lambda^1(M,\C)
\end{equation}
given by the Hodge connection $\Theta:\Lambda^1(X,\C) \to
\rho^*\Lambda^1(M,\C)$. This splitting induces a bigrading on the de
Rham algebra $\Lambda^\hdot(X,\C)$. The complex
$\Lambda^\hdot(X,\C)$ with this bigrading is a bicomplex which we
will denote by $\Lambda^{\hdot,\hdot}_\Theta(X)$. More precisely,
the de Rham differential $d$ is a sum of two anti-commuting
differentials
$$
\wt{d}:\Lambda^{\hdot,\hdot}_\Theta(X) \to
\Lambda^{\hdot+1,\hdot}_\Theta(X) \qquad\qquad
D:\Lambda^{\hdot,\hdot}_\Theta(X) \to
\Lambda^{\hdot,\hdot+1}_\Theta(X)
$$ 
(this is essentially equivalent to the flatness of the Hodge
connection $\Theta$).

Since the first term in the splitting \eqref{biTheta} is
$\Lambda^{1,0}_J(X)$, the subcomplexes $\Lambda^{\geq
k,\hdot}_\Theta(X)$ and $\Lambda^{\geq k,\hdot}_J(X)$ of the de Rham
complex $\Lambda^\hdot(X,\C)$ are the same for every $k$. Therefore
the associated graded quotients $\Lambda^{k,\hdot}_\Theta(X)$ and
$\Lambda^{k,\hdot}_J(X)$ are also isomorphic for every $k$. In
particular, we have
$$
\CC^{\hdot+2}(M) = \rho_*\Lambda^{2,0}_J(X) \cong
\rho_*\Lambda^{2,\hdot}_\Theta(X).
$$
On the other hand, since
$$
\Lambda^{1,0}_\Theta(X) = \Lambda^1(X,\C)/\rho^*\Lambda^1(M,\C) \cong
\Lambda^1(X/M,\C)
$$
is the bundle of relative $1$-forms on $X/M$, the quotient complex
$\Lambda^{\hdot,0}_\Theta(X)$ with the differential $\wt{d}$ is
canonically isomorphic
\begin{equation}\label{thtspl}
\Lambda^{\hdot,0}_\Theta(X) \cong \Lambda^\hdot(X/M,\C)
\end{equation}
to the relative de Rham complex $\Lambda^\hdot(X/M,\C)$. We use this
identification and choose as
$$
L^k(M) \subset \rho_*\Lambda^\hdot(X/M,\C) \cong
\rho_*\Lambda^{\hdot,0}_\Theta(X) 
$$
the subbundle of $k$-forms which are constant along the fibers of
the projection $\rho:X \to M$ (by a constant $k$-form on a vector
space $V$ with a basis $e_1,\dots,e_n$ we mean a linear combination
of forms $e_{a_1} \wedge \dots \wedge e_{a_k}$ with constant
coefficients). 

This choice guarantees that the relative de Rham differential
$\wt{d}:L^\hdot(M) \otimes \B^0(M) \to L^{\hdot+1}(M) \otimes
\B^0(M)$ takes a very simple form. Namely, it vanishes on $L^k(M)$,
and induces an isomorphism
\begin{equation}\label{SL}
\wt{d}:S^1(M) \cong L^1(M)
\end{equation}
between the generator subbundles $S^1(M) \subset \B^0(M)$ and
$L^1(M) \subset L^1(M) \otimes \B^1(M)$. This is important because
by construction $\wt{d}$ anti-commutes with $D$. Moreover, since
$\wt{d}$ obviously preserves the augmentation degrees, it
anti-commutes separately with each of the components $D_k$. Since we
already know the derivations $D_k$ on $\B^\hdot(M)$, the isomorphism
\eqref{SL} will allow us to compute individual components
$D_k:L^\hdot(M) \to L^\hdot(M) \otimes \B^1_k(M)$.

The first result is the following: the map $D_1 = D^{1,1}_1$ in
\eqref{d1} is minus the connection
$$
\nabla:L^{1,1}(M) \to L^{1,1}(M) \otimes \Lambda^1(M,\C)
$$
on the bundle $L^{1,1}(M) \cong \Lambda^{1,1}(M,\C)$. Indeed,
$D_1:L^\hdot(M) \to L^\hdot(M) \otimes \Lambda^1(M,\C)$ is a
derivation of the exterior algebra $L^\hdot(M)$. Thus it suffices to
prove that $D_1 = -\nabla$ on $L^1(M)$. Since $\wt{d}D_1 =
-D_1\wt{d}$, this follows from \eqref{SL} and the construction of
the map $D:\B^\hdot(M) \to \B^{\hdot+1}(M)$ given in
Section~\ref{pf}. This shows that taking
\begin{equation}\label{omega0}
\Omega_0 = \omega \in L^{1,1}(M) \cong \Lambda^{1,1}(M)
\end{equation}
guarantees that $D_1\Omega_0 = 0$.

At this point we will also choose the isomorphism $L^k(M) \cong
\Lambda^1(M,\C)$ -- namely, we take the composition of the embedding
$L^k(M) \subset \Lambda^k(X/M,\C) \cong \Lambda^{k,0}_\Theta(X)$ and
the restriction $i^*\Lambda^{k,0}(X) \to \Lambda^k(M,\C)$ to the
zero section $i:M \hookrightarrow X$. Then the form $\Omega_0$
defined by \eqref{omega0} automatically restricts to $\omega$.

It remains to prove that $D_2^l\Omega_0 = 0$. This is a corollary of
the following claim.

\begin{lemma}
The map $D_2^l$ defined in \eqref{d2l} is the composition of the
curvature
$$
R:L^{1,1}(M) \to L^{1,1}(M) \otimes \Lambda^{1,1}(M)
$$
of the connection $\nabla = D_1:L^{1,1}(M) \to L^{1,1}(M) \otimes
\Lambda^1(M,\C)$ and a certain bundle map
$$
Q:L^{1,1}(M) \otimes \Lambda^{1,1}(M) \to L^{0,2}(M) \otimes
\B^{2,-1}_{2,0}(M).
$$
\end{lemma}

\proof{} Extend the map $D_2^l$ to an algebra derivation
$$
D_2^l:L^\hdot(M) \to L^\hdot(M) \otimes \B^{2,-1}_{2,0}(M) \cong
L^\hdot(M) \otimes S^{1,-1}(M) \otimes \Lambda^{1,0}(M)
$$
by setting $D_2^l = 0$ on $L^{0,1}(M)$ and taking as
$$
D_2^l:L^{1,0}(M) \to L^{0,1}(M) \otimes \B^{2,-1}_{2,0}(M)
$$
the corresponding component of the map $D_2:L^{1,0}(M) \to L^1(M)
\otimes \B^1_2(M)$. Since $D_2^l$ is a derivation of the algebra
$L^\hdot(M)$ which vanishes on $L^{0,1}(M)$, on $L^{1,1}(M)$ it is
equal to the composition
$$
\begin{CD}
L^{1,0}(M) \otimes L^{0,1}(M) @>{D_2^l \otimes \id}>>
\B^{2,-1}_{2,0}(M) \otimes L^{0,1}(M) \otimes L^{0,1}(M) @>>>\\
@>{\id \otimes \Alt}>> \B^{2,-1}_{2,0}(M) \otimes L^{0,2}(M)
\end{CD}
$$
(here $\Alt:L^{0,1}(M) \otimes L^{0,1}(M) \to \Lambda^{0,2}(M)$ is
the alternation map). Therefore it suffices to prove that on
$L^{1,0}(M)$ we have
$$
D_2^l = P \circ R:L^{1,0}(M) \to L^{1,0}(M) \otimes
\Lambda^{1,1}(M) \to L^{0,1} \otimes \B^{2,-1}_{2,0}(M)
$$
for a certain bundle map $P:L^{1,0}(M) \otimes \Lambda^{1,1}(M) \to
L^{0,1} \otimes \B^{2,-1}_{2,0}(M)$. Since $D_2\wt{D} = -\wt{d}D_2$,
this follows directly from \eqref{SL} and \eqref{even.k} with
$k=2$. \endproof

This Lemma implies that $D_2^l(\Omega_0) = Q(R(\omega)) =
Q(\nabla(\nabla(\omega))) = 0$. This finishes the proof of
Theorem~\ref{hk}.

The last application of the formalism that we have developed in
this Section is the proof of Lemma~\ref{obata}.

\proof[Proof of Lemma~\ref{obata}.] The Obata connection $\nabla_O$
on a hypercomplex manifold $X$ is defined as follows: for every
$(0,1)$-form $\alpha \in \Lambda^{0,1}(X)$ the $(1,0)$-part
$\nabla_O^{1,0}$ is equal to
$$
\nabla_O^{1,0}\alpha = \6\alpha \in \Lambda^{0,1}(M) \otimes
\Lambda^{1,0}(M) ,
$$
while the $(0,1)$-part satisfies
$$
\nabla_O^{0,1}\alpha = -(\id \otimes j)(\6(j(\alpha))) \in
\Lambda^{0,1}(M) \otimes \Lambda^{0,1}(M).
$$
The first condition in fact automatically follows from the absence
of torsion.

Consider a hypercomplex structure on $X = \overline{T}M$
corresponding to a torsion-free connection $\nabla$ on $M$, and let
$\alpha \in \Lambda^{0,1}(M)$ be a $(0,1)$-form on $M$. We have to
prove that $\nabla_O^{0,1}(\rho^*\alpha) = \rho^*\nabla^{0,1}\alpha$
on the zero section $i:M \hookrightarrow X$. It suffices to prove
that
\begin{equation}\label{ob}
i^*\bar\6j(\rho^*\alpha) = i^*(\id \otimes
j)\left(\nabla^{0,1}\alpha\right)
\end{equation}
as section of the bundle $\Lambda^{0,1}(M) \otimes
i^*\Lambda^1(X,\C)$ on the zero section $M \subset X$. The canonical
isomorphism $H:\Lambda^{\hdot,\hdot}(X) \cong
\Lambda^{\hdot,\hdot}_J(X)$ sends the Dolbeault differential
$\bar\6$ to the component $D^{1,0}$ of the Dolbeault differential $D
= \bar\6_J:\Lambda^{\hdot,\hdot}_J(X) \to
\Lambda^{\hdot,\hdot+1}_J(X)$. Moreover, the composition of the map
$H \circ j:\rho^*\Lambda^{0,1}(M) \to \Lambda^{1,0}_J(X)$ and the
restriction to the zero section $M \subset X$ induces a bundle map
$$
P:\Lambda^{0,1}(M) \to L^1(M) = i^*\Lambda^{1,0}_J(X) \subset
i^*\Lambda^1(X,\C).
$$
It is easy to check that this map is proportional to the canonical
embedding $\Lambda^{0,1}(M) \cong L^{0,1}(M) \hookrightarrow
L^1(M)$. Therefore it commutes with the connection $\nabla^{1,0}$ --
namely, we have $(\id \otimes P) \circ \nabla^{1,0} = \nabla^{1,0}
\circ P$. The equation \eqref{ob} becomes
$$
D^{1,0}_1P(\alpha) = \nabla^{1,0}P(\alpha).
$$
But we have already proved that on $L^1(M)$ we have $\nabla =
D_1$. \endproof

\section{Symmetric spaces.}\label{symmm}

To illustrate our rather abstract methods by a concrete example, we
would like now to derive a formula for the canonical hypercomplex
structure on $X=\overline{T}M$ in the case when $M$ is a symmetric
space. In this case the classic equality $\nabla R = 0$ bring many
simplifications, so that the constructions of Section~\ref{pf} can
be seen through to a reasonably explicit final result. The formula
that we obtain is similar to the one obtained by O. Biquard and
P. Gauduchon in \cite{BG}.

Let us introduce some notations. Let $M$ be a symmetric space with
Levi-Civita connection $\nabla$ and curvature $R$. Consider the
total space $X =\overline{T}M$ with the canonical projection
$\rho:\overline{T}M \to M$. Let $A:\rho^*\Lambda^1(M) \to
\rho^*\Lambda^1(M)$ be the endomorphism of the pullback bundle
$\rho^*\Lambda^1(M)$ given by
\begin{equation}\label{A}
A(\xi)(\alpha) = -\frac{1}{3}R_m(\alpha) \cntrct (\xi \otimes \xi),
\qquad\qquad \alpha \in \rho^*\Lambda^1(M).
\end{equation}
Here $(\xi,m)$, $m \in M$, $\xi \in T_mM$ is a point in
$\overline{T}M$, $R_m$ is the curvature evaluated at the point $m
\in M$, and $R(\alpha)$ is interpreted as a section of the bundle
$\rho^*\Lambda^1(M) \otimes \rho^*\Lambda^{1,1}(M)$. Let also
$$
f(z) = \sum_{p \geq 1}f_pz^p
$$
be the generating function for the recurrence relation
\begin{equation}\label{rec}
f_p = -\frac{1}{2p+1}\sum_{1 \leq l \leq p-1}f_lf_{p-l}
\end{equation}
with the initial condition $f_1 = 1$. In other words, $f(z)$ is the
solution of the ODE
$$
2zf'(z) + f(z) + f^2(z) = 3z
$$
with the initial condition $f(0) = 0$. With these notations, we can
formulate our result.

\begin{prop}\label{symm}
Let $I:\Lambda^1(M) \to \Lambda^1(M)$ be the complex structure
operator on $M$. Then the map $J:\Lambda^1(X) \to \Lambda^1(X)$ for
the canonical normalized hypercomplex structure on $X$ is given by
the matrix
$$
J = \begin{pmatrix} 0 & f(A)I \\ I(f(A))^{-1} & 0 \end{pmatrix}
$$
with respect to the decomposition 
$$
\Lambda^1(X) = \Lambda^1(X/M) \oplus \rho^*\Lambda^1(M) \cong
\rho^*\Lambda^1(M) \oplus \rho^*\Lambda^1(M)
$$
associated to the Levi-Civita connection on $M$.
\end{prop}

We have already noted that this result is very similar to the
formula obtained in \cite{BG}. However, it is not the same. The
reason is the following: Biquard and Gauduchon work with the
cotangent bundle $T^*M$, and they use the normalization natural to
the cotangent bundle. As the result, their analog of the map $A$ is
slightly different, and the function $f(z)$ is also different -- in
particular, it is given by an explicit expression. To compare our
results with those of \cite{BG}, one should either compute the
normalization map $\LL:T^*M \to \overline{T}M$ for the
Biquard-Gauduchon hypercomplex structure, or go the other way around
and compute the map $\eta:\overline{T}M \to T^*M$ provided by
Lemma~\ref{cotang}. Unfortunately, I have not been able to do
either.

\proof[Proof of Proposition~\ref{symm}.] Throughout the proof, we
will freely use the notation of preceding Sections.

The main simplification in the case of symmetric spaces is the
well-known equality $\nabla R = 0$. For our construction it
immediately implies that the odd augmentation degree components
$D_{2p+1}$ of the weakly Hodge derivation $D:\B^\hdot(M) \to
\B^{\hdot+1}(M)$ vanish for $p \geq 1$. The only non-trivial
component is $D_1$. By \eqref{odd}, this immediately shows that the
connection $\Theta_{Re}$ on the fibration $\rho:X \to M$ is simply
the linear connection associated to the Levi-Civita connection
$\nabla$. Applying Lemma~\ref{j.matrix}, we see that all we have to
prove is the equality $R_J = f(A)$. By Lemma~\ref{Rj} this can be
rewritten in terms of the Weil algebra $\B^\hdot(M)$ as
$$
D_+\circ\sigma = f(A):\Lambda^1(M,\C) \to \B^1(M).
$$
Moreover, by \eqref{even.k} with $k=2$ the endomorphism $A$ of the
$\B^0(M)$-module $\B^1(M) \cong \rho_*\rho^*\Lambda^1(M,\C)$
satisfies
$$
A = D_2 \circ \sigma
$$
on $\Lambda^1(M,\C)$. In fact, this should be taken as the
definition -- I apologize to the reader for any possible mistakes in
writing down the explicit formula \eqref{A}.

Next we note that since $D_2$ vanishes on $\Lambda^1(M,\C)$, the
endomorphism $A$ is in fact equal
$$
A = \{D_2,\sigma\}
$$
to the commutator $\{D_2,\sigma\}$. Moreover, this formula holds not
only on the generator subbundle $\Lambda^1(M,\C) \subset \B^1(M)$,
but on the whole $\B^0(M)$-module $\B^1(M)$.  Indeed, by the
normalization condition \eqref{norm.bis} this commutator vanishes on
the generator subbundle $S^1(M) \subset \B^0(M)$. Since it is a
derivation of the Weil algebra, it vanishes on the whole
$\B^0(M)$. Therefore it restricts to a map of $\B^0(M)$-modules on
$\B^1(M) \subset \B^1(M)$.

We now trace one-by-one the induction steps in the proof of
Proposition~\ref{hc.prop}. Since all the odd terms vanish, we only
need to consider the even terms $D_{2p}$. We have to prove that
$$
\{D_{2p}, \sigma\} = f_kA^k:\B^1(M) \to \B^1(M).
$$
Since both sides are maps of $\B^0(M)$-modules, it suffices to prove
this on $\Lambda^1(M,\C) \subset \B^1(M)$, where we can replace
$\{D_{2p},\sigma\}$ with $D_p \circ \sigma$. By \eqref{even.k}, on
$S^1(M) \subset \B^0(M)$ we have
$$
D_{2p} = -\frac{1}{2p+1}\sum_{1 \leq l \leq p-1} \sigma \circ D_{2l}
\circ D_{2(p-l)}.
$$
Comparing this to \eqref{rec}, we see that it suffices to prove that
$$
\sigma \circ D_m \circ D_n \circ \sigma = \{D_m, \sigma\} \circ
\{D_n, \sigma\}:\Lambda^1(M) \to \B^1(M)
$$
for all even $m,n \geq 2$. Writing out the commutators on the
right-hand side, we see that the difference is equal to
$$
D_m \circ \sigma \circ D_n \circ \sigma + D_m \circ \sigma \circ
\sigma \circ D_n + \sigma \circ D_m \circ \sigma \circ D_n.
$$
The first summand vanishes since by the normalization condition
\eqref{norm.bis} we have $\sigma \circ D_n = 0$ on $S^1(M) =
\sigma(\Lambda^1(M,\C))$. The last two summands vanish since $D_n =
0$ on $\Lambda^1(M,\C)$. \endproof

\end{document}